\documentclass[11pt]{article}%
\usepackage{amsmath}
\usepackage{amsfonts}
\usepackage{amssymb}
\usepackage{graphicx}%
\usepackage{tikz}

\providecommand{\U}[1]{\protect\rule{.1in}{.1in}}
\providecommand{\U}[1]{\protect\rule{.1in}{.1in}}
\providecommand{\U}[1]{\protect\rule{.1in}{.1in}}
\RequirePackage[colorlinks,citecolor=blue,urlcolor=blue]{hyperref}
\providecommand{\U}[1]{\protect\rule{.1in}{.1in}}
\numberwithin{equation}{section}
\definecolor{ProcessBlue}{cmyk}{1,0,0,0.40}

\newtheorem{lemma}{Lemma}[section]
\newtheorem{theorem}{Theorem}[section]
\newtheorem{corollary}{Corollary}[section]
\newtheorem{proposition}{Proposition}[section]
\newtheorem{definition}{Definition}[section]

\newtheorem{remark}{Remark}[section]
\begin{document}

\title{\vspace{-1in}\parbox{\linewidth}{\footnotesize\noindent
} \vspace{\bigskipamount} \\Differential games with asymmetric information \\and without Isaacs' condition\thanks{The work is partially supported by the
Commission of the European Communities under the project SADCO,
FP7-PEOPLE-2010-ITN, No 264735, the French National Research Agency
ANR-10-BLAN 0112 and Natural Science Foundation of Jiangsu Province and China
(No.BK20140299; No.14KJB110022; No.11401414) and the collaborative innovation
center for quantitative calculation and control of financial risk.} }
\author{Rainer Buckdahn$^{1,3}$, Marc Quincampoix$^{1}$, Catherine Rainer$^{1}$,
Yuhong Xu$^{2,1}$\\$^{1}$Laboratoire de Math\'{e}matiques, CNRS-UMR 6205, Universit\'{e} de
Brest, France.\\$^{2}$Mathematical center for interdiscipline research and School of
mathematical sciences,\\Soochow University, Suzhou 215006, P. R. China.\\$^{3}$School of Mathematics, Shandong University, Jinan 250100, P. R. China.\\E-mails: rainer.buckdahn@univ-brest.fr; marc.quincampoix@univ-brest.fr; \\catherine.rainer@univ-brest.fr; yuhong.xu@hotmail.com}
\date{Feb. 16th, 2015}
\maketitle

\indent\textbf{Abstract.} We investigate a two-player zero-sum differential
game with asymmetric information on the payoff and without Isaacs' condition.
The dynamics is an ordinary differential equation parametrized by two
controls chosen by the players. Each player has a private information on the
payoff of the game, while his opponent knows only the probability distribution
on the information of the other player.

We show that a suitable definition of random strategies allows to prove the
existence of a value in mixed strategies.  This value is taken in the sense of the limit of any time discretization, as the mesh of the time partition tends to zero. We characterize it
 in terms of the unique viscosity solution in some
dual sense of a Hamilton-Jacobi-Isaacs equation. Here we do not suppose the
Isaacs' condition, which is usually assumed in differential games.

\vspace{0.3in}\indent\textbf{Key words.} zero-sum differential game;
asymmetric information; Isaacs' condition; viscosity solution; subdynamic
programming principle; dual game \newline


\section{Introduction}

We consider a two-player differential game which dynamic is given by the
following differential equation%

\[
\frac{dX_{s}}{ds } =f\left(  X_{s},u_{s},v_{s}\right)  \text{, } \;
s\in\left[  t,T\right]
\]
parametrised by two controls $u : \left[  t,T\right]  \mapsto U $ and $v :
\left[  t,T\right]  \mapsto V$ chosen by the Players 1 and 2, respectively.
Here the function $f : \mathbb{R}^{d} \times U \times V \mapsto\mathbb{R} ^{d}
$ satisfies standard assumptions, and $U$ and $V$ denote two compact metric spaces.

The final cost involves $I \times J$ payoffs $g_{ij}(X_{T})$ $i=1,2 \ldots I$,
$j=1,2 \ldots J$, where the $g_{ij}:\mathbf{R}^{d}\mapsto\mathbf{R}$ are given
bounded functions. The first player aims to minimize the cost, while the
second player's objective is to maximize it.

Let us now describe how the game is played. For this purpose let us fix an
initial time $t\in\lbrack0,T]$.

\noindent- Before the game starts, a pair $(i,j)$ is chosen randomly according
to a probability measure $p\otimes q\in\Delta(I)\times\Delta(J)$. Here
$\Delta(I)$ denotes the set of probabilities $p=(p_{i})_{i=1,\ldots,I}$ on
\{$1,\ldots,I$\}; $\Delta(J)$ is similarly defined.

\noindent- The choice of $i$ is communicated to Player 1 but not to Player 2,
while $j$ is communicated to Player II but not to Player I;

\noindent- The game is played on the time interval $[t,T]$;

\noindent- Both players know the probability $p\otimes q$ and observe their
opponents' controls during the game.

\noindent Note that the players do not know which $g_{ij}$ they are actually
optimising, because they have only a part of the information on the pair
$(i,j)$. Nevertheless they can try to guess their missing information by
observing what their opponent does. Indeed, in order to use his information, a
player necessarily reveals at least a part of it, and any piece of information
he reveals can be later exploited by his opponent.

\medskip

Games with asymmetric information where introduced by Aumann and Maschler in the 1960s \textrm{\cite{am}}. They proved that games with lack of information on one side have a uniform value. The existence of a value in the general case but in a weaker sense was established by Mertens and Zamir \textrm{\cite{mz}}. 
Models with asymmetric information in economics have been
investigated extensively, for instance in \textrm{\cite{b93, bbhps, beb, zie}%
}. They study the case that participants in a market have private information
not public to the others.

\medskip

A milestone in the literature of differential games is the article
\textrm{\cite{es} }which has been later extended to stochastic differential
games in \textrm{\cite{fs}}. It was shown there that under the following
Isaacs' condition
\[
\underset{u\in U}{\inf}\underset{v\in V}{\sup}f\left(  x,u,v\right)  \cdot
\xi=\underset{v\in V}{\sup}\underset{u\in U}{\inf}f\left(  x,u,v\right)
\cdot\xi,
\]
the value function of a differential game is given as the unique viscosity
solution of a Hamilton--Jacobi--Isaacs (HJI in short) equation. For two-person
zero-sum stochastic differential games, we also refer the reader to
\textrm{\cite{bcq}} for an overview and a more complete description.
Differential games with asymmetric information on the payoff were studied
first by Cardaliaguet in \textrm{\cite{c07}} ( see also \textrm{\cite{cr091,
cr12}}) . The case where the asymmetric information concerns the initial
position was studied in \cite{CJQ}. The extension to stochastic differential
games was investigated in \textrm{\cite{cr092, cr13}}. The proof was
accomplished by introducing the notion of dual viscosity solutions to the HJI
equation of a usual differential game, where the probability $\left(
p,q\right)  $ just appear as additional parameters. Such a notion of dual
solution was introduced in \cite{c07} for differential games and in \cite{BdM}
for repeated games. A different unique characterisation via the viscosity
solution of the HJI equation with double obstacles in the form of constraints
in $(p,q)$ was given in \textrm{\cite{c09}}. More recently, Oliu-Barton
\textrm{\cite{oliubarton} }extends Cardaliaguet \cite{c07} to the case of
correlated types, and where the lack of information carries over the initial
state, the dynamic and the pay-off function.

Differential games without Isaacs' assumption where first considered by Krasovskii and Subbotin\textrm{\cite{K}}, using relaxed controls. This approach has also been chosen by Sirbu \cite{sirbu}. 
Recently, the article \textrm{\cite{blq1} } considered zero-sum differential
games with complete information without Isaacs' condition by playing classical, randomized controls, imposing on the
underlying controls for both players a conditional independence property. The article
\textrm{\cite{blq2}} generalized this method to stochastic differential
games using the approach of backward stochastic differential equations.\\
In comparison with \cite{blq1} we have made here the choice of a randomisation of the strategies applying pathwise to control processes, rather than to consider strategies applying globally (and not only pathwise) to randomised control processes. Both approaches lead, of course, to the same value, which also can be gotten by considering relaxed controls. 
While the approach of \cite{blq1} with its randomisation of controls is in some sense still near to this idea of relaxed controls, the choice we have made in the present work is justified by the fact that both players use strategies, the associated controls are only a consequence and computed with the help of the chosen strategies. 
In addition, our approach is not only justified from the point of view of interpretation but it leads also to a considerable simplification of the computation. 

\medskip

The present paper introduces this method to differential games with asymmetric
information. By using a suitable notion of random, non-anticipative strategies
with delay, we show that the upper and lower value
function of the game along a partition of $[t,T]$ (where the players are restricted
to react at times $t=t_{0}<t_{1}<\ldots<t_{N}=T$) converge to the unique dual
viscosity solution of the following HJI equation
\begin{equation}
\left\{
\begin{array}
[c]{ll}%
\frac{\partial V}{\partial t}\left(  t,x\right)  + H\left(  x,DV\left(  t,x\right)\right)  =0\text{,}
& \text{in }\left[  0,T\right]  \times\mathbf{R}^{d}\text{,}\\
V\left(  T,x\right)  =\sum_{ij}p_{i}q_{j}g_{ij}\left(  x\right)  \text{,} &
\end{array}
\right.  \label{1.1}%
\end{equation}
where%
\[
H\left(  x,\xi\right)  =
\underset{\mu\in
\mathcal{P}\left(  U\right)  }{\inf}\underset{\nu\in\mathcal{P}\left(
V\right)  }{\sup}\int_{U\times V}
f\left(  x,u,v\right)  \mu\left(  du\right)  \nu\left(  dv\right)  \cdot\xi \; .
\]
Here $\mathcal{P}\left(  U\right)  $ denotes the space of all probability
measures on $U$, $\mathcal{P}\left(  V\right)  $ that on $V$.  Since both control state spaces $U$ and $V$ are compact, ${\cal P}(U)$ and  ${\cal P}(V)$ are compact, convex subsets of the linear topological spaces of signed measures over $U$ and $V$, respectively. Moreover, taking into acount the continuity of $f(x,.,.):U\times V\rightarrow \mathbf{R}$, for all $(\mu,\nu)\in{\cal P}(U)\times{\cal P}(V)$ the functions 
${\cal P}(U)\ni\mu'\rightarrow \int_{U\times V}f(x,u,v)\mu'(du)\nu(dv)$ and ${\cal P}(V)\ni\nu'\rightarrow \int_{U\times V}f(x,u,v)\mu(du)\nu'(dv)$
are continuous (w.r.t. the topology generated by the weak convergence of measures), and from the linearity of the functions we get their convexity and concavity, respectively. But then it follows, for instance, from Sion's Minimax Theorem, that 
\begin{equation}
\label{isaacs}
\underset{\mu\in\mathcal{P}\left(  U\right)  }{\inf}\underset{\nu
\in\mathcal{P}\left(  V\right)  }{\sup}\int_{U\times V}
f\left(  x,u,v\right)  \mu\left(  du\right)  \nu\left(  dv\right)  \cdot\xi.
 =\underset{\nu\in\mathcal{P}\left(  V\right)  }{\sup
}\underset{\mu\in\mathcal{P}\left(  U\right)  }{\inf}\int_{U\times V}
f\left(  x,u,v\right)  \mu\left(  du\right)  \nu\left(  dv\right)  \cdot\xi \;
 .
\end{equation}

\bigskip

Non-anticipative strategies are randomised on a single probability space. The
obtained limit value is similar to the value with relaxed controls concerning
$\mathcal{P}\left(  U\right)  $ and $\mathcal{P}\left(  V\right)  $ as the
control spaces, although players are not allowed to use relaxed controls. See
\textrm{\cite{blq1}} for more explanation.

The organisation of the paper is as follows. In the next section we present
necessary definitions. Section 3 shows the convexity (concavity) of the value
of the game along a partition and considers its dual game. A subdynamic
programming principle is established for the dual game in section 4. In
section 5, we prove that the lower\ value function and the upper value
function are a viscosity subsolution and a viscosity supersolution
respectively of the associated HJI equation. In section 6, we show by a
special comparison principle of partial differential equations that the limit
value of the game along partitions exists as the mesh of partitions tends to
zero, and the limit value function is characterized as the dual solution of
some HJI equation. In the last section we consider the case of lack of
information on the dynamics.

\section{Settings of the game}

Throughout the paper we work with the following probability space:%
\[
(\Omega,\mathcal{F},P):=\left(  \left[  0,1\right]  ,\mathcal{B}\left(
\left[  0,1\right]  \right)  ,dx\right)  ,
\]
where the interval $[0,1]$ is endowed with the Borel $\sigma$-field
$\mathcal{B}\left(  \left[  0,1\right]  \right)  $ and the Borel measure $dx.$
Let $\left\{  \zeta_{j,l},\, l\ge1,\ j=1,2\right\}  $ be a family of
independent random variables following all a uniform distribution on $[0,1].$
Let us consider two compact metric spaces $U$ and $V$ representing the control
state spaces used by player 1 and 2, respectively. $\mathcal{P}\left(
U\right)  $ and $\mathcal{P}\left(  V\right)  $ denote the space of all
probability measures over $U$ and $V$, endowed with Borel $\sigma$-field
$\mathcal{B}\left(  U\right)  $ and $\mathcal{B}\left(  V\right)  $,
respectively. It is an immediate consequence of Skorohod's Representation
Theorem that $\mathcal{P}\left(  U\right)  \ $(resp., $\mathcal{P}\left(
V\right)  $) coincides with the set of the laws of all $U$-valued (resp.,
$V$-valued) random variables defined over $\left(  \left[  0,1\right]
,\mathcal{B}\left(  \left[  0,1\right]  \right)  ,dx\right)  $.

Let us now introduce the admissible controls for both players.

\begin{definition}
(Admissible controls) Given the initial time $t\in\left[  0,T\right]  $, we
define the sets of admissible controls for player 1 and player 2 by

$\mathcal{U}_{t,T}=\left\{  \text{all U-valued and Lebesque measurable
functions }\left(  u_{s}\right)  _{s\in\left[  t,T\right]  }\right\}  $,

$\mathcal{V}_{t,T}=\left\{  \text{all V-valued and Lebesque measurable
functions}\left(  v_{s}\right)  _{s\in\left[  t,T\right]  }\right\}  $.

\noindent Both spaces are endowed with the topology generated by the
convergence in $L^{1}([t,T]).$
\end{definition}

The dynamics of the game is given by%
\begin{equation}
\left\{
\begin{array}
[c]{ll}%
dX_{s}^{t,x,u,v}=f\left(  X_{s}^{t,x,u,v},u_{s},v_{s}\right)  ds\text{,} &
s\in\left[  t,T\right]  \text{, }\left(  u,v\right)  \in\mathcal{U}%
_{t,T}\times\mathcal{V}_{t,T},\\
X_{t}^{t,x,u,v}=x\text{,} &
\end{array}
\right.  \label{2.1}%
\end{equation}
where $f:\mathbf{R}^{d}\times U\times V\mapsto\mathbf{R}^{d}$ is bounded,
continuous, and Lipschitz continuous in $x\in\mathbf{R}^{d}$, uniformly with
respect to $u$ and $v$. Standard estimates show that there exists a constant
$C>0$ such that, for all $\left(  t,x\right)  ,\left(  t^{\prime},x^{\prime
}\right)  \in\left[  0,T\right]  \times\mathbf{R}^{d}$ and all $s\in\left[
t\vee t^{\prime},T\right]  $,

(i) $\left\vert X_{s}^{t,x,u,v}-x\right\vert \leq C\left(  s-t\right)  $,

(ii) $\left\vert X_{s}^{t,x,u,v}-X_{s}^{t^{\prime},x^{\prime},u,v}\right\vert
\leq C\left(  \left\vert t-t^{\prime}\right\vert +\left\vert x-x^{\prime
}\right\vert \right)  $.

While player 1 tries to minimise a given cost functional, the objective of
player 2 is to maximise it. For this they have at their disposal the above
introduced spaces of admissible controls. But can they play the game ``control
against control''? As it is by now well-known, apart from rather particular
cases, differential games of the type ``control against control'' don't, in
general, admit the dynamic programming principle and don't have a value. This
is why in the literature approaches studying games of the type ``strategy
against control'' and ``non anticipative strategy with delay (NAD-strategy)
against NAD-strategy'' have imposed. This latter type has turned out to be the
best adapted one for the study of differential games with asymmetric
information. Taking into account the asymmetry of the information, the players
aim to hide a part of their private knowledge. To do this they randomise their
strategies. From a technical point of view, this randomness is also the key
argument to get a value of our game without Isaacs' condition.

Let us consider now a partition $\pi=\left\{  0=t_{0}<t_{1}<\ldots
<t_{N}=T\right\}  $ and let us fix arbitrarily the initial time of the game
$t\in\left[  t_{k-1},t_{k}\right]  $, for some $0\leq k\leq N$.

\begin{definition}
\label{def2.2}(Random NAD strategies along the partition $\pi$) A random
NAD-strategy along the partition $\pi$ for player 1 is a mapping
$\alpha:\Omega\times\mathcal{V}_{t,T}\mapsto\mathcal{U}_{t,T}$ of the form
\[
\alpha\left(  \omega,v\right)  \left(  s\right)  =\alpha_{l}\left(  \left(
\zeta_{1,k},\ldots,\zeta_{1,l}\right)  \left(  \omega\right)  ,v\right)
\left(  s\right)  ,
\]
$\omega\in\Omega$, $s\in\lbrack t\vee t_{l-1},t\vee t_{l} )$, $k\leq l\leq N$,
with a Borel measurable mapping $\alpha_{l}:\mathbf{R}^{l-k+1}\times
\mathcal{V}_{t,T}\mapsto\mathcal{U}_{t,T}$ satisfying the following: For all
$v,v^{\prime}\in\mathcal{V}_{t,T},$ it holds that, whenever $v=v^{\prime}$
a.e. on $\left[  t,t_{l-1}\right]  $, we have for all $x\in\mathbf{R}^{l-k+1}%
$, $\alpha_{l}\left(  x,v\right)  =\alpha_{l}\left(  x,v^{\prime}\right)  $,
a.e on $\left[  t\vee t_{l-1},t\vee t_{l}\right]  $.

Similarly, a random NAD strategy along the partition $\pi$ for player 2 is a
mapping $\beta:\Omega\times\mathcal{U}_{t,T}\mapsto\mathcal{V}_{t,T}$ of the
form
\[
\beta\left(  \omega,u\right)  \left(  s\right)  =\beta_{l}\left(  \left(
\zeta_{2,k},\ldots,\zeta_{2,l}\right)  \left(  \omega\right)  ,u\right)
\left(  s\right)  ,
\]
$\omega\in\Omega$, $s\in\left[  t\vee t_{l-1},t\vee t_{l}\right]  $, $k\leq
l\leq N$, with a Borel measurable mapping $\beta_{l}:\mathbf{R}^{l-k+1}%
\times\mathcal{U}_{t,T}\mapsto\mathcal{V}_{t,T}$ satisfying the following: For
all $u,u^{\prime}\in\mathcal{U}_{t,T},$ it holds that, whenever $u=u^{\prime}$
a.e. on $\left[  t,t_{l-1}\right]  $, we have for all $x\in\mathbf{R}^{l-k+1}%
$, $\beta_{l}\left(  x,u\right)  =\beta_{l}\left(  x,u^{\prime}\right)  $, a.e
on $\left[  t\vee t_{l-1},t\vee t_{l}\right]  $.
\end{definition}

We denote by $\mathcal{A}_{r}^{\pi}\left(  t,T\right)  $ the set of all such
random NAD strategies for player 1, and by $\mathcal{B}_{r}^{\pi}\left(
t,T\right)  $ that for player 2. Sometimes we will use pure (i.e.
deterministic) strategies: Let $\mathcal{A}^{\pi}\left(  t,T\right)  $ (resp.
$\mathcal{B}^{\pi}\left(  t,T\right)  $) denote the subset of strategies in
$\mathcal{A}_{r}^{\pi}\left(  t,T\right)  $ (resp. $\mathcal{B}_{r}^{\pi
}\left(  t,T\right)  $) which do not depend on $\omega\in\Omega$.

\begin{remark}
\label{rem2.1}Given $t\in\left[  t_{k-1},t_{k}\right]  $, for $s\in\left[
t,T\right]  $, $\left(  u,v\right)  \in\mathcal{U}_{t,T}\times\mathcal{V}%
_{t,T}$, we have
\begin{equation}%
\begin{array}
[c]{c}%
\alpha\left(  \omega,v\right)  \left(  s\right)  =\alpha_{k}\left(
\zeta_{1,k}\left(  \omega\right)  ,v\right)  \cdot\mathbf{1}_{[t,t_{k}%
)}\left(  s\right)  +\sum_{l=k+1}^{N}\alpha_{l}\left(  \left(  \zeta
_{1,k},\ldots,\zeta_{1,l}\right)  \left(  \omega\right)  ,v\right)
\cdot\mathbf{1}_{[t_{l-1},t_{l})}\left(  s\right)  ,\\
\beta\left(  \omega,u\right)  \left(  s\right)  =\beta_{k}\left(  \zeta
_{2,k}\left(  \omega\right)  ,u\right)  \cdot\mathbf{1}_{[t,t_{k})}\left(
s\right)  +\sum_{l=k+1}^{N}\beta_{,l}\left(  \left(  \zeta_{2,k},\ldots
,\zeta_{2,l}\right)  \left(  \omega\right)  ,u\right)  \cdot\mathbf{1}%
_{[t_{l-1},t_{l})}\left(  s\right)  .
\end{array}
\label{alphabeta}%
\end{equation}

\end{remark}

\begin{remark}
The above notion of strategies comes close to the one used in the framework of repeated games. One difference is that here the actions of the players on each interval  $[t_l,t_{l+1}]$ depend on the whole sequence of random variables $\zeta_{i,k},\ldots,\zeta_{i,l}$ introduced since the beginning of the game, while in repeated games, the strategy at stage $l$ depends only on the last, new random variable. \\
In fact, in the proof of the dynamic programming principles -crucial step to establish the HJI-equation-, it is important that the strategies may depend on the past trajectory of the dynamic. There are several possibilities to take this into account: firstly,  as for repeated games, by letting the action at stage $l$ depend not only on the action of the opponent but also on its own past actions, or, secondly, by observing directly this trajectory. In our setting, the most convenient was to keep, at each stage $l$ the memory of the randomness introduced in the previous steps. In terms of the required information, all these approaches are equivalent.
\end{remark}

Thanks to the delay of the strategies, we have the following property.

\begin{lemma}
\label{lem2.1}For any $(\alpha,\beta)\in\mathcal{A}_{r}^{\pi}\left(
t,T\right)  \times\mathcal{B}_{r}^{\pi}\left(  t,T\right)  $ there is a unique
(up to a null set) mapping $\Omega\ni\omega\mapsto\left(  u_{\omega}%
,v_{\omega}\right)  \in\mathcal{U}_{t,T}\times\mathcal{V}_{t,T}$ such that,
for all $\omega\in\Omega,$
\begin{equation}
\label{fixpoint}
\alpha\left(  \omega,v_{\omega}\right)  =u_{\omega}\text{, }\beta\left(
\omega,u_{\omega}\right)  =v_{\omega}\text{, a.e. on }\left[  t,T\right]
\text{.}%
\end{equation}

\end{lemma}

The proof of this lemma uses standard arguments. However since the special
form of the random NAD-strategies which we use here is new, we prefer to give
the proof for the convenience of the reader.

\paragraph{Proof}

Let $t\in\left[  t_{k-1},t_{k}\right]  $, $0\leq k\leq N$. For each $\omega
\in\Omega,\, \alpha\left(  \omega,v\right)  $ and $\beta\left(  \omega
,u\right)  $ restricted to $\left[  t,t_{k}\right]  $ depend only on
$v\in\mathcal{V}_{t,T}$ and $u\in\mathcal{U}_{t,T}$ restricted to $\left[
t,t_{k-1}\right]  $. But $\left[  t,t_{k-1}\right]  $ is empty or a singleton,
so that $\alpha\left(  \omega,v\right)  $, $\beta\left(  \omega,u\right)  $ on
$\left[  t,t_{k}\right]  $ do not depend on $v$ and $u$. Thus, for any
$v^{0}\in\mathcal{V}_{t,T}$ and $u^{0}\in\mathcal{U}_{t,T}$, we can define
$u_{\omega}^{1}=\alpha\left(  \omega,v^{0}\right)  $, $v_{\omega}^{1}%
=\beta\left(  \omega,u^{0}\right)  $. We observe that this definition
guarantees the measurability of the mapping $\Omega\ni\omega\mapsto\left(
u_{\omega}^{1},v_{\omega}^{1}\right)  \in\mathcal{U}_{t,T}\times
\mathcal{V}_{t,T}.$ Moreover, from this definition we have
\[
\alpha\left(  \omega,v^{1}\right)  =u^{1}\text{, }\beta\left(  \omega
,u^{1}\right)  =v^{1}\text{, a.e., on }\left[  t,t_{k}\right]  \text{.}%
\]

Assume that for $j\geq2$, $\Omega\ni\omega\mapsto\left(  u_{\omega}^{j-1},
v_{\omega}^{j-1}\right)  \in\mathcal{U}_{t,T}\times\mathcal{V}_{t,T}$ is
measurable and such that $\alpha\left(  \omega,v_{\omega}^{j-1}\right)
=u_{\omega}^{j-1}$, $\beta\left(  \omega,u_{\omega}^{j-1}\right)  =v_{\omega
}^{j-1}$, a.e., on $\left[  t,t_{k+j-2}\right]  $. Then we set $u_{\omega}%
^{j}=\alpha\left(  \omega,v_{\omega}^{j-1}\right)  $, $v_{\omega}^{j}%
=\beta\left(  \omega,u_{\omega}^{j-1}\right)  $. Obviously, the thus defined
mapping $\omega\mapsto\left(  u_{\omega}^{j},v_{\omega}^{j}\right)  $ is
measurable and $\left(  u_{\omega}^{j},v_{\omega}^{j}\right)  =\left(
u_{\omega}^{j-1},v_{\omega}^{j-1}\right)  $, a.e., on $\left[  t,t_{k+j-2}%
\right]  $. Then by the NAD property of $\alpha$ and $\beta$, $u_{\omega}%
^{j}=\alpha\left(  \omega,v_{\omega}^{j}\right)  $, $v_{\omega}^{j}%
=\beta\left(  \omega,u_{\omega}^{j}\right)  $, a.e., on $\left[
t,t_{k+j-1}\right]  $. Iterating the above steps, we can obtain the desired
result. The uniqueness is an immediate consequence of the above construction.
$\Box$\\

 Thanks to Lemma \ref{lem2.1}, to any pair $(\alpha
,\beta)\in\mathcal{A}_{r}^{\pi}(0,T)\times\mathcal{B}_{r}^{\pi}(0,T)$ and any
$(x,\omega)\in\mathbf{R}^{d}\times\lbrack0,1]$ can be associated a trajectory
$t\mapsto X_{t}^{t,x,\alpha(\omega,\cdot),\beta(\omega,\cdot)}$ defined as:
\[
X_{t}^{t,x,\alpha(\omega,\cdot),\beta(\omega,\cdot)}:=X_{t}^{t,x,u_{\omega
},v_{\omega}},
\]
with $(u_{\omega},v_{\omega})$ uniquely determined by relation (\ref{fixpoint}).

\begin{remark}
We observe that, for all $1\le l\le N-1$, the processes $u$ and $v$
constructed in the above proof and restricted to the time interval $[t,t_{l}]$
are conditionally independent knowing $\zeta_{k}=(\zeta_{1,k},\zeta
_{2,k}),\dots,\zeta_{l-1}=(\zeta_{1,l-1},\zeta_{2,l-1})$. Indeed, the
processes $u$ and $v$ are of the following form:
\[
\left\{
\begin{array}
[c]{l}%
u_{\cdot}\left(  s\right)  =u^{k}\left(  s, \zeta_{1,k}\right)  \cdot
\mathbf{1}_{[t,t_{k})}\left(  s\right)  +\sum_{l=k+1}^{n}u^{l}\left(  s,
\zeta_{k},\ldots,\zeta_{l-1},\zeta_{1,l}\right)  \cdot\mathbf{1}%
_{[t_{l-1},t_{l})}\left(  s\right)  \text{,}\\
v_{\cdot}\left(  s\right)  =v^{k}\left(  s, \zeta_{2,k}\right)  \cdot
\mathbf{1}_{[t,t_{k})}\left(  s\right)  +\sum_{l=k+1}^{n}v^{l}\left(  s,
\zeta_{k},\ldots,\zeta_{l-1},\zeta_{2,l}\right)  \cdot\mathbf{1}%
_{[t_{l-1},t_{l})}\left(  s\right)  \text{,}%
\end{array}
\right.
\]
\bigskip where $\left(  u^{l},v^{l}\right)  $ are measurable functions of $s$,
$\zeta_{j,m}$, $j=1,2$, $k\leq m\leq l-1$ and $\zeta_{l}=\left(  \zeta
_{1,l},\zeta_{2,l}\right)  $, for $l\geq1$.
\end{remark}

The description of the game involves $I\times J$ terminal payoffs (where
$I,J\geq1$): $g_{ij}:\mathbf{R}^{d}\mapsto\mathbf{R}$ for $i=1,\ldots,I$ and
$j=1,\ldots,J$, which are supposed to be Lipschitz continuous and bounded
throughout the paper.\newline

We now define the lower and the upper value functions.

\begin{definition}
Let $\left(  p,q\right)  \in\Delta\left(  I\right)  \times\Delta\left(
J\right)  $, $\left(  t,x\right)  \in\lbrack0,T)\times\mathbf{R}^{d}$. For
$t\in\lbrack t_{k-1},t_{k})$, $\hat{\alpha}=\left(  \alpha_{i}\right)
_{i=1,\ldots I}\in\left(  \mathcal{A}_{r}^{\pi}\left(  t,T\right)  \right)
^{I}$, $\hat{\beta}=\left(  \beta_{j}\right)  _{j=1,\ldots J}\in\left(
\mathcal{B}_{r}^{\pi}\left(  t,T\right)  \right)  ^{J}$, we define the cost
functional
\begin{equation}
\mathcal{J}\left(  t,x,\hat{\alpha},\hat{\beta},p,q\right)  =\sum_{i=1}%
^{I}\sum_{j=1}^{J}p_{i}q_{j}E\left[  g_{ij}\left(  X_{T}^{t,x,\alpha_{i}%
,\beta_{j}}\right)  \right]  \text{, }%
\end{equation}
the upper value function%
\begin{equation}
W^{\pi}\left(  t,x,p,q\right)  =\underset{\hat\alpha\in\left(  \mathcal{A}%
_{r}^{\pi}\left(  t,T\right)  \right)  ^{I}}{\inf}\underset{\hat\beta
\in\left(  \mathcal{B}_{r}^{\pi}\left(  t,T\right)  \right)  ^{J}}{\sup
}\mathcal{J}\left(  t,x,\hat{\alpha},\hat{\beta},p,q\right)  \text{, }%
\end{equation}
as well as the lower value function%
\begin{equation}
V^{\pi}\left(  t,x,p,q\right)  =\underset{\hat\beta\in\left(  \mathcal{B}%
_{r}^{\pi}\left(  t,T\right)  \right)  ^{J}}{\sup}\underset{\hat\alpha
\in\left(  \mathcal{A}_{r}^{\pi}\left(  t,T\right)  \right)  ^{I}}{\inf
}\mathcal{J}\left(  t,x,\hat{\alpha},\hat{\beta},p,q\right)  \text{. }%
\end{equation}

\end{definition}

\section{Convexity and Fenchel conjugates}

The convexity (resp. concavity) of the value functions with respect to the
probabilities $p$ (resp. to $q$) is a crucial aspect in the analysis of games
with asymmetric information. In this section we study this property and
introduce the Fenchel conjugate of the value functions.\newline

The proof of the following lemmas are inspired by \textrm{\cite{c07}}.
We just remark that, in order to take into account the missing Isaacs'
assumption, our notion of strategy is much more explicit and, therefore, more
restrictive than the one defined in \textrm{\cite{c07}}. This explains
additional difficulties in the proofs.

\begin{lemma}
\label{lem3.1} The value functions $W^{\pi}$ and $V^{\pi}$ are Lipschitz
continuous with respect to $\left(  t,x,p,q\right)  $, uniformly with respect
to $\pi$.
\end{lemma}

\paragraph{Proof}

We only give the proof for $V^{\pi}$. By the definition of $V^{\pi}$ and the
boundness of $g_{ij}$, it is easy to see that $V^{\pi}$ is Lipschitz with
respect $p$ and $q$. For every $t\in\lbrack0,T]$ and $\left(  u,v\right)
\in\mathcal{U}_{t,T}\times\mathcal{V}_{t,T}$, we can show that $x\rightarrow
g_{ij}\left(  X_{T}^{t,x,u,v}\right)  $ is Lipschitz uniformly w.r.t.
$(t,u,v)$ and, hence, for any $\left(  \hat{\alpha},\hat{\beta}\right)
\in\left(  \mathcal{A}_{r}^{\pi}\left(  t,T\right)  \right)  ^{I}\times\left(
\mathcal{B}_{r}^{\pi}\left(  t,T\right)  \right)  ^{J}$, the mapping%
\[
x\mapsto\mathcal{J}\left(  t,x,\hat{\alpha},\hat{\beta},p,q\right)
\]
is Lipschitz continuous with a Lipschitz constant $C$ (for short,
$C$-Lipschitz) independent of $\left(  t,p,q\right)  \in\lbrack0,T]\times
\Delta\left(  I\right)  \times\Delta\left(  J\right)  $ and $\pi$. Therefore,
we deduce that $V^{\pi}$ is $C$-Lipschitz with respect to $x$.

\bigskip Now we show that $V^{\pi}$ is Lipschitz in $t$. Let $x\in
\mathbf{R}^{d}$, $\left(  p,q\right)  \in\Delta\left(  I\right)  \times
\Delta\left(  J\right)  $ and $t<t^{\prime}<T$ be fixed. Let $\hat{\beta
}=\left(  \beta_{j}\right)  _{j=1,\ldots J}\in\left(  \mathcal{B}_{r}^{\pi
}\left(  t,T\right)  \right)  ^{J}$ be $\varepsilon$-optimal for $V^{\pi
}\left(  t,x,p,q\right)  $. We have to associate with each $\beta_{j}$ a
strategy $\beta_{j}^{\prime}\in\mathcal{B}_{r}^{\pi}\left(  t^{\prime
},T\right)  $. To this aim, we fix some arbitrary constant control $\bar{u}\in
U$ and set
\[
\tilde{\beta}_{j}\left(  \omega,u\right)  =\beta_{j}\left(  \omega,\tilde
{u}\right)  \text{, where }\tilde{u}\left(  s\right)  =\left\{
\begin{array}
[c]{ll}%
\bar{u}\text{,} & s\in\lbrack t,t^{\prime})\text{,}\\
u(s)\text{,} & s\in\left[  t^{\prime},T\right]  \text{.}%
\end{array}
\right.
\]
If $t^{\prime}<t_{k}$, then $\tilde{\beta}_{j}\in\mathcal{B}_{r}^{\pi}\left(
t^{\prime},T\right)  $ and we set $\beta_{j}^{\prime}=\tilde{\beta}_{j}%
$.\newline Otherwise, let $l>k$ such that $t_{l-1}\leq t^{\prime}<t_{l}$. We
consider now $l-k+1$ random variables $\kappa_{k},\ldots,\kappa_{l}$ on
$([0,1],\mathcal{B}([0,1]),dx)$ with $\kappa_{l}(x)=x,x\in\lbrack0,1]$,
which are uniformly distributed on $[0,1]$, mutually
independent, and independent of $\zeta_{i,m},(i,m)\neq(2,l)$. We remark that
then also the composed random variables $\kappa_{k}\circ\zeta_{2,l}%
,\ldots,\kappa_{l}\circ\zeta_{2,l}$ are mutually independent, uniformly
distributed random variables, which are moreover independent of all
$\zeta_{i,m},(i,m)\neq(2,l)$. Indeed, recall that the $\zeta_{i,l}$'s
themselves are also i.i.d. and obey a uniform distribution over the interval
$[0,1]$.\newline Now we set, for any $u\in\mathcal{U}_{t^{\prime},T}$, $l\leq
m\leq N$, and $s\in\lbrack t^{\prime}\vee t_{m-1},t_{m})$,
\[
\beta_{j}^{\prime}(\omega,u)(s)=\tilde{\beta}_{j,m}\left(  \kappa_{k}%
\circ\zeta_{2,l},\ldots,\kappa_{l}\circ\zeta_{2,l},\zeta_{2,l+1},\ldots
,\zeta_{2,m}\right)  (\omega,u)(s).
\]
Then $\beta_{j}^{\prime}\in\mathcal{B}_{r}^{\pi}\left(  t^{\prime},T\right)
$. Moreover, for all $u\in\mathcal{U}_{t^{\prime},T}$,$\beta_{j}^{\prime}(u)$
and $\tilde{\beta}_{j}(u)$ obey the same law knowing $\zeta_{2,l}%
,\zeta_{2,l+1},\ldots,\zeta_{2,n}$, if $t_{l-1}\leq t^{\prime}<t_{l}$ (Here we use the fact that $\kappa_l(\zeta_{2,l})=\zeta_{2,l}$). Thus
\[
E[g(X_{T}^{t^{\prime},x,u,\beta_{j}^{\prime}(u)})]=E[g(X_{T}^{t^{\prime
},x,u,\tilde{\beta}_{j}(u)})],
\]
and, hence, for all $\hat{\alpha}=(\alpha_{i})\in\left(  \mathcal{A}_{r}^{\pi
}\left(  t^{\prime},T\right)  \right)  ^{I}$,
\[
\mathcal{J}(t^{\prime},x,\hat{\alpha},(\beta_{j}^{\prime}),p,q)=\mathcal{J}%
(t^{\prime},x,\hat{\alpha},(\tilde{\beta}_{j}),p,q).
\]
Next, to every $\alpha\in\mathcal{A}_{r}^{\pi}\left(  t^{\prime},T\right)  $,
we associate some new strategy $\alpha^{\prime}\in\mathcal{A}_{r}^{\pi}\left(
t,T\right)  $ by setting, for all $v\in\mathcal{V}_{t,T}$,
\[
\alpha^{\prime}\left(  \omega,v\right)  \left(  s\right)  =\left\{
\begin{array}
[c]{ll}%
\bar{u}\left(  s\right)  \text{,} & s\in\lbrack t,t^{\prime})\text{,}\\
\alpha\left(  \omega,v|_{\left[  t^{\prime},T\right]  }\right)  \left(
s\right)  \text{,} & s\in\left[  t^{\prime},T\right]  \text{.}%
\end{array}
\right.
\]
By construction, if $t_{l-1}\leq t^{\prime}\leq t_{l}\;(k\leq l\leq n)$, the
couples of random controls associated by Lemma \ref{lem2.1} to the couples of
strategies $(\alpha^{\prime},\beta_{j})$ and $(\alpha,\tilde{\beta}_{j})$
coincide on the time interval $[t^{\prime},T]$ under the conditional law
$P[\;\cdot\;|\zeta_{2,l},\zeta_{l+1},\ldots,\zeta_{n}]$ . Therefore the
standard estimate applies :
\[
E\left[  |X_{s}^{t,x,\alpha^{\prime},\beta_{j}}-X_{s}^{t^{\prime}%
,x,\alpha,\tilde{\beta}_{j}}|\right]  \leq M\left\vert t^{\prime}-t\right\vert
,s\in\lbrack t^{\prime},T],
\]
where $M$ is a constant depending on the bound of $f$ as well as the Lipschitz
constants of $f$ and the functions $g_{ij}$, but not on $\pi$. Hence, for any
$\hat{\alpha}=\left(  \alpha_{i}\right)  _{i=1,\ldots I}\in\left(
\mathcal{A}_{r}^{\pi}\left(  t^{\prime},T\right)  \right)  ^{I}$, we have%
\begin{align*}
\mathcal{J}\left(  t^{\prime},x,\hat{\alpha},\left(  \beta_{j}^{\prime
}\right)  ,p,q\right)  =\mathcal{J}\left(  t^{\prime},x,\hat{\alpha},\left(
\tilde{\beta}_{j}\right)  ,p,q\right)   &  \geq\mathcal{J}\left(
t,x,\hat{\alpha}^{\prime},\hat{\beta},p,q\right)  -LM\left\vert t^{\prime
}-t\right\vert \\
&  \geq\underset{\hat{\alpha}^{\prime\prime}\in\left(  \mathcal{A}_{r}^{\pi
}\left(  t,T\right)  \right)  ^{I}}{\inf}\mathcal{J}\left(  t,x,\hat{\alpha
}^{\prime\prime},\hat{\beta},p,q\right)  -LM\left\vert t^{\prime}-t\right\vert
\\
&  \geq V^{\pi}\left(  t,x,p,q\right)  -\varepsilon-LM\left\vert t^{\prime
}-t\right\vert
\end{align*}
where we have used the fact that $\hat{\beta}$ is $\varepsilon$-optimal for
$V^{\pi}\left(  t,x,p,q\right)  $. Consequently,%
\begin{equation}
V^{\pi}\left(  t^{\prime},x,p,q\right)  \geq V^{\pi}\left(  t,x,p,q\right)
-\varepsilon-LM\left\vert t^{\prime}-t\right\vert . \label{3.1}%
\end{equation}
To prove the reverse inequality, we associate in a symmetric way to above to a
vector of strategies $\hat{\beta}=\left(  \beta_{j}\right)  _{j=1,\ldots J}%
\in\left(  \mathcal{B}_{r}^{\pi}\left(  t^{\prime},T\right)  \right)  ^{J}$
which is $\varepsilon$-optimal for $V^{\pi}(t^{\prime},x,p,q)$ and each
arbitrary $\hat{\alpha}\in\left(  \mathcal{A}_{r}^{\pi}\left(  t,T\right)
\right)  ^{I}$ some vectors of strategies $\hat{\beta}^{\prime}\in\left(
\mathcal{B}_{r}^{\pi}\left(  t,T\right)  \right)  ^{J}$ and $\hat{\alpha
}^{\prime}\in\left(  \mathcal{A}_{r}^{\pi}\left(  t^{\prime},T\right)
\right)  ^{I}$, in order to get the inequality
\[
\mathcal{J}\left(  t,x,\hat{\alpha},\hat{\beta}^{\prime},p,q\right)
\geq\mathcal{J}\left(  t^{\prime},x,\hat{\alpha}^{\prime},\hat{\beta
},p,q\right)  -LM|t-t^{\prime}|\geq V^{\pi}\left(  t^{\prime},x,p,q\right)
-\varepsilon-LM\left\vert t^{\prime}-t\right\vert
\]
and, thus,
\begin{equation}
V^{\pi}\left(  t,x,p,q\right)  \geq V^{\pi}\left(  t^{\prime},x,p,q\right)
-\varepsilon-LM\left\vert t^{\prime}-t\right\vert . \label{3.2}%
\end{equation}

Finally, thanks to the arbitrariness of $\varepsilon>0$, from \eqref{3.1} and
\eqref{3.2} we get the Lipschitz continuity for $V^{\pi}$ in $t$. $\Box$\\

 Next we have to prove that $V^\pi$ and $W^\pi$ are convex in $p$ and concave in $q$.
We can use the standard argumentation of \cite{sorin}, provide we show that the sets of strategies satisfy the following convexity property.

\begin{lemma}
Let $(t,x,p)\in[0,T)\times\mathbf{R}^d\times\Delta(I)$, $q^0,q^1\in\Delta(J)$ and $\lambda\in (0,1)$. Set $q^\lambda=(1-\lambda)q^0+\lambda q^1$. For all $\hat\beta^0,\hat\beta^1\in
\left(  \mathcal{B}_{r}^{\pi}\left(  t,T\right)  \right)^{J}$, 
there exists $\hat\beta^\lambda\in\left(  \mathcal{B}_{r}^{\pi}\left(  t,T\right)  \right)  ^{J}$, 
such that, for all $\hat\alpha\in\left(\mathcal{A}_{r}^{\pi}\left(t,T\right)\right)^I$,
\[ \mathcal{J}\left(  t,x,\hat\alpha,\hat{\beta}^{\lambda},p,q^{\lambda
}\right) =(1-\lambda)\mathcal{J}(  t,x,\hat\alpha,\hat{\beta}^0,p,q^0) +\lambda\mathcal{J}(  t,x,\hat\alpha,\hat{\beta}^1,p,q^1).
\]
An analogue result holds for the strategies of Player 1.
\end{lemma}

\paragraph{Proof}

 For all $j\in\{ 1,\ldots, J\}$, we set  $c_{j}%
=\frac{\left(  1-\lambda\right)  q_{j}^{0}}{q_{j}^{\lambda}}$.
For $\omega\in\Omega$, $u\in\mathcal{U}_{t,T}$, $s\in\lbrack t\vee
t_{l-1},t\vee t_{l})$, $k\leq l\leq N$, $j=1,\ldots,J$, we define the
following strategies%
\[
\beta_{j}^{\lambda}\left(  \omega,u\right)  \left(  s\right)  =\beta
_{lj}\left(  \left(  \zeta_{2,k},\ldots,\zeta_{2,l}\right)  \left(
\omega\right)  ,u\right)  \left(  s\right)
\]
where
\begin{align*}
\beta_{lj}\left(  \left(  y_{1},\ldots,y_{l-k+1}\right)  ,u\right)   &
=\beta_{lj}^{0}\left(  \frac{1}{c_{j}}y_{1},y_{2},\ldots,y_{l-k+1},u\right)
\cdot\mathbf{1}_{\left[  0,c_{j}\right]  }\left(  y_{1}\right) \\
&  +\beta_{lj}^{1}\left(  \frac{1}{1-c_{j}}\left(  y_{1}-c_{j}\right)
,y_{2},\ldots,y_{l-k+1},u\right)  \cdot\mathbf{1}_{\left[  c_{j},1\right]
}\left(  y_{1}\right)  .
\end{align*}
The mappings $\beta_{l,j}^{0}$ and $\beta_{l,j}^{1}$, $k\leq l\leq N$, are
associated with $\beta_{j}^{0}$ and $\beta_{j}^{1}$ through Definition 2.2. It
is easy to see that $\left(  \beta_{j}^{\lambda}\right)  \in\left(
\mathcal{B}_{r}^{\pi}\left(  t,T\right)  \right)  ^{J}$.\newline
Fix now $(i,j)\in\{ 1,\ldots,I\}\times\{ 1,\ldots,J\}$. A straight forward computation yields
\[\begin{array}{l}
E[g_{i,j}(X^{t,x,\alpha_i,\beta^\lambda_j})]\\
=\displaystyle\int_{[0,c_j]}
E\left[ g_{i,j}\left(X^{t,x,\alpha_i,\beta^0_j\left((\frac 1{c_j}y_1,\zeta_{2,k+1},\ldots,\zeta_{2,N})(\omega)\right)}\right)\right]dy_1\\
\ \ \ \
+\displaystyle\int_{[c_j,1]}E\left[ g_{i,j}\left(X^{t,x,\alpha_i,\beta^1_j\left((\frac 1{1-c_j}(y_1-c_j),\zeta_{2,k+1},\ldots,\zeta_{2,N})(\omega)\right)}\right)\right]dy_1\\
=\frac{\left(
1-\lambda\right)  q_{j}^{0}}{q_{j}^{\lambda}}E\left[  g_{ij}\left(
X_{T}^{t,x,\alpha_i,\beta_{j}^{0}}\right)  \right]  +\frac{\lambda q_{j}^{1}%
}{q_{j}^{\lambda}}E\left[  g_{ij}\left(  X_{T}^{t,x,\alpha_i,\beta_{j}^{1}%
}\right)  \right]
.
\end{array}\]
Consequently it holds that
\[\begin{array}{l}
\sum_{j=1}^Jq^\lambda_jE\left[  g_{ij}\left(  X_{T}^{t,x,\alpha_i,\beta_{j}^{\lambda}}\right)  \right]
\\
=
\left(  1-\lambda\right)  \sum_{j=1}^{J}q_{j}^{0}E\left[  g_{ij}\left(  X_{T}^{t,x,\alpha_i,\beta_{j}^{0}}\right)  \right] 
+\lambda\sum_{j=1}^{J}q_{j}^{1}E\left[  g_{ij}\left(
X_{T}^{t,x,\alpha_i,\beta_{j}^{1}}\right)  \right]
\end{array}\]
The result follows.
 $\Box$\newline 

\begin{lemma}  (see \cite{sorin}, Chapter 2) 
\label{lem3.2}For any $\left(  t,x\right)  \in\lbrack0,T)\times\mathbf{R}^{d}%
$, the mappings $W^{\pi}\left(  t,x,p,q\right)  $, $V^{\pi}\left(
t,x,p,q\right)  $ are convex in $p$ and concave in $q$ on $\Delta\left(
I\right)  $ and $\Delta\left(  J\right)  $, respectively.
\end{lemma}

\noindent For the proof the reader is referred to \cite{sorin}.\\

Let us now introduce the Fenchel conjugates: Given a mapping $w:[0,T)\times
\mathbf{R}^{d}\times\Delta(I)\times\Delta(J)\mapsto\mathbf{R}$ convex in $p$
and concave in $q$ on $\Delta\left(  I\right)  $ and $\Delta(J)$,
respectively, we denote by $w^{\ast}$ its convex conjugate with respect to
variable $p$:%
\[
w^{\ast}\left(  t,x,\hat{p},q\right)  :=\underset{p\in\Delta\left(  I\right)
}{\sup}\left\{  \left\langle \hat{p},p\right\rangle -w\left(  t,x,p,q\right)
\right\}  \text{, }\left(  t,x,\hat{p},q\right)  \in\lbrack0,T]\times
\mathbf{R}^{d}\times\mathbf{R}^{I}\times\Delta\left(  J\right)  ,
\]
and $w^{\#}$ its concave conjugate with respect to variable $q$:%
\[
w^{\#}\left(  t,x,p,\hat{q}\right)  :=\underset{q\in\Delta\left(  J\right)
}{\inf}\left\{  \left\langle \hat{q},q\right\rangle -w\left(  t,x,p,q\right)
\right\}  \text{, }\left(  t,x,p,\hat{q}\right)  \in\lbrack0,T]\times
\mathbf{R}^{d}\times\Delta\left(  I\right)  \times\mathbf{R}^{J}.
\]
In coherence with these notations, we write $V^{\pi\ast}$ ($W^{\pi\#})$ for
the convex (resp. concave) conjugate of $V^{\pi}$ ($W^{\pi})$ with respect to
$p$ (resp. $q$).  As the previous lemma, the following alternative formulation of $V^{\pi\ast}$ is standard for normal-form games with convex sets of strategies (see also \cite{sorin}, Chapter 2). For this reason we omit the proof.

\begin{lemma}
\label{lem3.3}(Reformulation of $V^{\pi\ast}$) For all $\left(  t,x,\hat
{p},q\right)  \in\lbrack0,T)\times\mathbf{R}^{d}\times\mathbf{R}^{I}%
\times\Delta\left(  J\right)  $,\newline(i) we have
\begin{equation}
 \label{3.5}
V^{\pi\ast}\left(  t,x,\hat{p},q\right)  =\underset{\left(  \beta_{j}\right)
\in\left(  \mathcal{B}_{r}^{\pi}\left(  t,T\right)  \right)  ^{J}}{\inf
}\underset{\alpha\in\mathcal{A}^{\pi}\left(  t,T\right)  }{\sup}\underset
{i\in\left\{  1,\ldots,I\right\}  }{\max}\left\{  \hat{p}_{i}-\sum_{j=1}%
^{J}q_{j}E\left[  g_{ij}\left(  X_{T}^{t,x,\alpha,\beta_{j}}\right)  \right]
\right\} 
\end{equation}

\end{lemma}

\section{The subdynamic programming principle}

 For any continuous time control problem, the first step leading to the Hamilton-Jacobi-equation is a dynamic programming principle. In the case of differential games with asymmetric information, we have to prove two sub- (resp. super-) dynamic programming principles: one for each Fenchel conjugate.\\

\begin{lemma}
\label{lem4.1}(Subdynamic programming principle) For any $\left(  t,x,\hat
{p},q\right)  \in\lbrack t_{k-1},t_{k})\times\mathbf{R}^{d}\times
\mathbf{R}^{I}\times\Delta\left(  J\right)  $ and for all l ($k\leq l\leq N$),
we have%
\begin{align}
V^{\pi\ast}\left(  t,x,\hat{p},q\right)   &  \leq\underset{\left(  \beta
_{j}\right)  \in\mathcal{B}_{r}^{\pi}\left(  t,t_{l}\right)  }{\inf}%
\underset{\alpha\in\mathcal{A}^{\pi}\left(  t,t_{l}\right)  }{\sup}E\left[
V^{\pi\ast}\left(  t_{l},X_{t_{l}}^{t,x,\alpha,\beta},\hat{p},q\right)
\right] \label{4.1}\\
&  \leq\underset{\left(  \beta_{j}\right)  \in\mathcal{B}_{r}^{\pi}\left(
t,t_{l}\right)  }{\inf}\underset{\alpha\in\mathcal{A}_{r}^{\pi}\left(
t,t_{l}\right)  }{\sup}E\left[  V^{\pi\ast}\left(  t_{l},X_{t_{l}}%
^{t,x,\alpha,\beta},\hat{p},q\right)  \right]  .
\end{align}

\end{lemma}

\paragraph{Proof}

The second inequality is obvious. Let us prove the first one. For this, we let
$V^{\pi\ast}(t,t_{l},x,\hat{p},q)$ denote the right side of \eqref{4.1}. For
arbitrarily given $\varepsilon>0$, let $\beta^{0}\in\mathcal{B}_{r}^{\pi
}\left(  t,t_{l}\right)  $ be an $\varepsilon$-optimal strategy for
$V^{\pi\ast}\left(  t,t_{l},x,\hat{p},q\right)  $. For any $z\in\mathbf{R}%
^{d}$, let $\hat{\beta}^{z}=\left(  \beta_{j}^{z}\right)  \in\left(
\mathcal{B}_{r}^{\pi}\left(  t_{l},T\right)  \right)  ^{J}$ be an
$\varepsilon$-optimal strategy for player 2 in the dual game with value
function $V^{\pi\ast}\left(  t_{l},z,\hat{p},q\right)  $. From the uniform
Lipschitz continuity of the mappings%
\[
y\mapsto\underset{\alpha\in\mathcal{A}_{r}^{\pi}\left(  t_{l},T\right)  }%
{\sup}\underset{i\in\left\{  1,\ldots,I\right\}  }{\max}\left\{  \hat{p}%
_{i}-\sum_{j=1}^{J}q_{j}E\left[  g_{ij}\left(  X_{T}^{t_{l},y,\alpha,\beta
_{j}^{x}}\right)  \right]  \right\}
\]
and, hence, also that of $y\mapsto V^{\pi\ast}\left(  t_{l},y,\hat
{p},q\right)  $, it follows that $\hat{\beta}^{z}$ is a $\left(
2\varepsilon\right)  $-optimal strategy for $V^{\pi\ast}\left(  t_{l}%
,y,\hat{p},q\right)  ,$ if $z\ $belongs to the ball $B_{r}\left(  y\right)  $,
for some radius $r>0$ small enough but not depending on $y\in\mathbf{R}^{d}$.
Because the coefficient $f$ is bounded, $X_{t_{l}}^{t,x,\alpha,\beta}$ takes
all its values in some ball $B_{R}\left(  0\right)  $, for $R>0$ large enough.
We choose a finite sequence $\left(  x_{n}\right)  _{n=1,\ldots,n_{0}}%
\subset\mathbf{R}^{d}$ such that $B_{R}\left(  0\right)  \subset\cup
_{n=1}^{n_{0}}B_{r}\left(  x_{n}\right)  $. This allows to construct a Borel
partition $\left(  A_{n}\right)  _{n=1,\ldots,n_{0}}$ of the ball
$B_{R}\left(  0\right)  $ such that, for all $1\leq n\leq n_{0}$, $x_{n}\in
A_{n}\subset B_{r}\left(  x_{n}\right)  $. To simplify the notation, we write
$\beta_{j}^{n}:=\beta_{j}^{x_{n}}$. We observe that on the event $\{X_{t_{l}%
}^{t,x,\alpha,\beta}\in B_{r}\left(  x_{n}\right)  \}$ the strategy $\beta
_{j}^{n}$ is $\left(  2\varepsilon\right)  $-optimal for $V^{\pi\ast}\left(
t_{l},X_{t_{l}}^{t,x,\alpha,\beta},\hat{p},q\right)  $.

For any $\omega\in\Omega$ and $u\in\mathcal{U}_{t,T}$, we set%
\[
\beta_{j}\left(  \omega,u\right)  \left(  s\right)  =\left\{
\begin{array}
[c]{ll}%
\beta^{0}\left(  \omega,u\right)  \left(  s\right)  \text{,} & s\in\lbrack
t,t_{l})\text{,}\\
\beta_{j}^{n}\left(  \omega,u|_{\left[  t_{l},T\right]  }\right)  \text{,} &
s\in\left[  t_{l},T\right]  \text{ and }X_{t_{l}}^{t,x,u,\beta^{0}}(\omega)\in
A_{n}\text{.}%
\end{array}
\right.
\]
Since, for all $1\leq n\leq n_{0},$ $(\omega,u)\mapsto\mathbf{1}_{A_{n}%
}\left(  X_{t_{l}}^{t,x,u,\beta^{0}}(\omega)\right)  $ is $\sigma\left\{
\zeta_{2,k},\ldots,\zeta_{2,l},u|_{[t,t_{l})}\right\}  $-measurable, there are
Borel functions $f^{n}:\mathbf{R}^{l-k+1}\times\mathcal{U}_{t,t_{l}}%
\mapsto\mathbf{R}$ such that $\mathbf{1}_{A_{n}}\left(  X_{t_{l}}%
^{t,x,u,\beta^{0}}\right)  =f^{n}\left(  \zeta_{2,k},\ldots,\zeta
_{2,l},u|_{[t,t_{l})}\right)  $. Note also that, by the definition of
$\beta_{j}^{n}\in\mathcal{B}_{r}^{\pi}\left(  t_{l},T\right)  $,
\[
\beta_{j}(\omega,u)(s)=\left(  \beta^{0}(\omega,u)(s)\mathbf{1}_{[t,t_{l}%
)}\left(  s\right)  +(\sum_{n=1}^{n_{0}}\beta_{j}^{n}(\omega,u|_{\left[
t_{l},T\right]  })(s)f^{n}\left(  \zeta_{2,k}(\omega),\ldots,\zeta
_{2,l}(\omega),u|_{[t,t_{l})}\right)  )\mathbf{1}_{\left[  t_{l},T\right]
}\left(  s\right)  \right)  .
\]
Thus $\beta_{j}\in\mathcal{B}_{r}^{\pi}\left(  t,T\right)  $.

For any $\alpha\in\mathcal{A}^{\pi}\left(  t,T\right)  $, we have%
\[
g_{ij}\left(  X_{T}^{t,x,\alpha,\beta_{j}}\right)  =\sum_{n=1}^{n_{0}}%
g_{ij}\left(  X_{T}^{t_{l},X_{t_{l}}^{t,x,\alpha,\beta^{0}},\tilde{\alpha
},\beta_{j}^{n}}\right)  \cdot\mathbf{1}_{A_{n}}\left(  X_{t_{l}}%
^{t,x,\alpha,\beta^{0}}\right)  ,
\]
where $\tilde{\alpha}\in$ $\mathcal{A}^{\pi}_r\left(  t_{l},T\right)  $ is a
restriction of $\alpha$ to $\left[  t_{l},T\right]  $ defined by%
\[
\tilde{\alpha}\left(  v\right)  \left(  s\right)  =\alpha\left(  v^{\prime
}\right)  \left(  s\right)  \text{, }\forall v\in\mathcal{V}_{t,T}\text{,
where }v^{\prime}\left(  s\right)  =\left\{
\begin{array}
[c]{ll}%
\bar{v}\left(  s\right)  \text{,} & s\in\lbrack t,t_{l})\text{,}\\
v\left(  s\right)  \text{,} & s\in\left[  t_{l},T\right]  \text{,}%
\end{array}
\right.
\]
with the controls $\left(  \bar{u}\left(  \cdot\right)  ,\bar{v}\left(
\cdot\right)  \right)  $ being associated with $\left(  \alpha,\beta
^{0}\right)  $ through Lemma \ref{lem2.1}.\\
 Remark that, since $\beta^0$ is a deterministic function of $(\zeta_{2,k},\ldots,\zeta_{l-1})$, the same holds for $\bar v$ and for $\tilde\alpha$. During the computations, we shall need to fix this dependence by the notation $\tilde\alpha:=\tilde\alpha(\zeta_{2,k},\ldots,\zeta_{l-1})$.\\
Further we observe that $X_{t_{l}}^{t,x,\alpha,\beta^{0}}$ is independent of
$\beta_{j}^{n}$. Indeed, while $\beta_{j}^{n}$ depends only on $(\zeta
_{2,l},\ldots,\zeta_{2,N})$, $X_{t_{l}}^{t,x,\alpha,\beta^{0}}=X_{t_{l}%
-}^{t,x,\alpha,\beta^{0}}$ only depends on $(\zeta_{2,k},\ldots,\zeta
_{2,l-1})$. It follows that
\[%
\begin{array}
[c]{l}%
\underset{i\in\left\{  1,\ldots,I\right\}  }{\max}\left\{  \hat{p}_{i}%
-\sum_{j=1}^{J}q_{j}E\left[  g_{ij}\left(  X_{T}^{t,x,\alpha,\beta_{j}%
}\right)  \right]  \right\} \\
=\underset{i\in\left\{  1,\ldots,I\right\}  }{\max}\left\{  \hat{p}_{i}%
-\sum_{j=1}^{J}q_{j}E\left[  \sum_{n=1}^{n_{0}}g_{ij}\left(  X_{T}%
^{t_{l},X_{t_{l}}^{t,x,\alpha,\beta^{0}},\tilde{\alpha},\beta_{j}^{n}}\right)
\cdot\mathbf{1}_{A_{n}}\left(  X_{t_{l}}^{t,x,\alpha,\beta^{0}}\right)
\right]  \right\} \\
=\underset{i\in\left\{  1,\ldots,I\right\}  }{\max}\left\{  \hat{p}_{i}%
-\sum_{j=1}^{J}q_{j} E\left[  \sum_{n=1}^{n_{0}}E\left[  g_{ij}\left(
X_{T}^{t_{l},y,\tilde{\alpha}(z),\beta_{j}^{n}}\right)  \right]  _{y=X_{t_{l}}^{t,x,\alpha,\beta^{0}},z=(\zeta_{2,k},\ldots,\zeta_{2,l-1})}\cdot\mathbf{1}_{A_{n}}\left(  X_{t_{l}}%
^{t,x,\alpha,\beta^{0}}\right)  \right]  \right\} \\
\leq E\left[  \sum_{n=1}^{n_{0}}\underset{i\in\left\{  1,\ldots,I\right\}
}{\max}\left\{  \hat{p}_{i}-\sum_{j=1}^{J}q_{j} E\left[  g_{ij}\left(
X_{T}^{t_{l},y,\tilde{\alpha}(z),\beta_{j}^{n}}\right)  \right]  \right\}
_{y=X_{t_{l}}^{t,x,\alpha,\beta^{0}},z=(\zeta_{2,k},\ldots,\zeta_{2,l-1})}\cdot\mathbf{1}_{A_{n}}\left(  X_{t_{l}%
}^{t,x,\alpha,\beta^{0}}\right)  \right] \\
\leq E\left[  \sum_{n=1}^{n_{0}}\underset{\alpha^{\prime}\in\mathcal{A}%
_{r}^{\pi}\left(  t_{l},T\right)  }{\sup}\underset{i\in\left\{  1,\ldots
,I\right\}  }{\max}\left\{  \hat{p}_{i}-\sum_{j=1}^{J}q_{j}E\left[
g_{ij}\left(  X_{T}^{t_{l},y,\alpha^{\prime},\beta_{j}^{n}}\right)  \right]
\right\}  _{y=X_{t_{l}}^{t,x,\alpha,\beta^{0}}}\cdot\mathbf{1}_{A_{n}}\left(
X_{t_{l}}^{t,x,\alpha,\beta^{0}}\right)  \right] \\
\leq E\left[  \sum_{n=1}^{n_{0}}V^{\pi\ast}\left(  t_{l},X_{t_{l}}%
^{t,x,\alpha,\beta^{0}},\hat{p},q\right)  \cdot\mathbf{1}_{A_{n}}\left(
X_{t_{l}}^{t,x,\alpha,\beta^{0}}\right)  \right]  +2\varepsilon\\
\leq V^{\pi\ast}\left(  t,t_{l},x,\hat{p},q\right)  +3\varepsilon.
\end{array}
\]
The latter inequality is due to the fact that $\beta^{0}$ is $\varepsilon
$-optimal for $V^{\pi\ast}\left(  t,t_{l},x,\hat{p},q\right)  $. Then we
conclude that $V^{\pi\ast}\left(  t,x,\hat{p},q\right)  \leq V^{\pi\ast
}\left(  t,t_{l},x,\hat{p},q\right)  $. $\Box$

\section{ Viscosity solutions of the dual game}

Let $\left(  \pi_{n}\right)  _{n\geq1}$ be a given sequence of partitions of
the time interval $\left[  0,T\right]  $ such that the mesh of the partition
$|\pi_{n}|$ tends to zero, when $n\rightarrow\infty$. We prove in this section
that $\left(  V^{\pi_{n}\ast}\right)  $ converges to some function $\tilde{V}$
as $\left\vert \pi_{n}\right\vert $ $\rightarrow0$ and that, for fixed
$\left(  \hat{p},q\right)  $, the limit function $\tilde{V}$ is a viscosity
subsolution of the following HJI equation:%
\begin{equation}
\frac{\partial\tilde{V}}{\partial t}\left(  t,x\right)  +H^{\ast}\left(
x,D\tilde{V}\left(  t,x\right)  \right)  =0\text{,} \text{in }\left[
0,T\right]  \times\mathbf{R}^{d}. \label{5.1}%
\end{equation}
where
\begin{align}
H^{\ast}\left(  x,\xi\right)   &  =-H\left(  x,-\xi\right) \nonumber\\
&  =\underset{\nu\in\mathcal{P}\left(  V\right)  }{\inf}\underset{\mu
\in\mathcal{P}\left(  U\right)  }{\sup}\left(  {\displaystyle\int
\nolimits_{U\times V}}f\left(  x,u,v\right)  \mu\left(  du\right)  \nu\left(
dv\right)  \cdot\xi\right) \nonumber\\
&  =\underset{\mu\in\mathcal{P}\left(  U\right)  }{\sup}\underset{\nu
\in\mathcal{P}\left(  V\right)  }{\inf}\left(  {\displaystyle\int
\nolimits_{U\times V}}f\left(  x,u,v\right)  \mu\left(  du\right)  \nu\left(
dv\right)  \cdot\xi\right)  . \label{5.2}%
\end{align}

The following lemma is obtained directly from the boundedness of $f$ and the
Lipschitz continuity of $f\left(  \cdot,u,v\right)  $ and of $g\left(
\cdot\right)  $.

\begin{lemma}
\label{lem5.1} There is some constant $L\in\mathbf{R}^{+}$, depending only on
the bound of f and the Lipschitz constant of $f\left(  \cdot,u,v\right)  $ and
of g, such that, for all partition $\pi$ of the interval $\left[  0,T\right]
$ and for all $\left(  t,x,\hat{p},q\right)  ,\left(  t^{\prime},x^{\prime
},\hat{p}^{\prime},q^{\prime}\right)  \in\lbrack0,T)\times\mathbf{R}^{d}%
\times\mathbf{R}^{I}\times\Delta\left(  J\right)  $,%
\[
\left\vert V^{\pi\ast}\left(  t,x,\hat{p},q\right)  -V^{\pi\ast}\left(
t^{\prime},x^{\prime},\hat{p}^{\prime},q^{\prime}\right)  \right\vert \leq
L\left(  \left\vert t-t^{\prime}\right\vert +\left\vert x-x^{\prime
}\right\vert +\left\vert \hat{p}-\hat{p}^{\prime}\right\vert +\left\vert
q-q^{\prime}\right\vert \right)  .
\]

\end{lemma}

By the above equi-Lipschitz continuity of the family of lower dual value
functions, applying the Arzel\`{a}-Ascoli Theorem, we have

\begin{lemma}
\label{lem5.2}There exists a subsequence of partitions, still denoted by
$\left(  \pi_{n}\right)  _{n\geq1}$ and there is a bounded Lipschitz functions
$\tilde{V}:\left[  0,T\right]  \times\mathbf{R}^{d}\times\mathbf{R}^{I}%
\times\Delta\left(  J\right)  \mapsto\mathbf{R}$ such that

$V^{\pi_{n}\ast}\rightarrow\tilde{V}$ uniformly on compact sets in $\left[
0,T\right]  \times\mathbf{R}^{d}\times\mathbf{R}^{I}\times\Delta\left(
J\right)  $.
\end{lemma}

We deduce from Lemma \ref{lem5.1} that also the function $\tilde{V}$ in Lemma
\ref{lem5.2} is Lipschitz continuous with respect to all its variables. More
precisely, we get the following corollary:

\begin{corollary}
\label{cor5.1}For the constant L introduced in Lemma \ref{lem5.1}, we have,\\
for all $\left(  t,x,\hat{p},q\right)  ,\left(  t^{\prime},x^{\prime},\hat
{p}^{\prime},q^{\prime}\right)  \in\lbrack0,T)\times\mathbf{R}^{d}%
\times\mathbf{R}^{I}\times\Delta\left(  J\right)  $,%
\[
\left\vert \tilde{V}\left(  t,x,\hat{p},q\right)  -\tilde{V}\left(  t^{\prime
},x^{\prime},\hat{p}^{\prime},q^{\prime}\right)  \right\vert \leq L\left(
\left\vert t-t^{\prime}\right\vert +\left\vert x-x^{\prime}\right\vert
+\left\vert \hat{p}-\hat{p}^{\prime}\right\vert +\left\vert q-q^{\prime
}\right\vert \right)  .
\]

\end{corollary}

 \begin{remark}
Lemma \ref{lem5.1} and Lemma \ref{lem5.2} also hold for $W^{\pi\#}$ on
$[0,T)\times\mathbf{R}^{d}\times\Delta\left(  I\right)  \times\mathbf{R}^{J}$.
Let $\tilde{W}:[0,T)\times\mathbf{R}^{d}\times\Delta\left(  I\right)
\times\mathbf{R}^{J}\mapsto\mathbf{R}$ denote the uniform limit of $W^{\pi
_{n}\#}$ on compact sets of $[0,T)\times\mathbf{R}^{d}\times\Delta\left(
I\right)  \times\mathbf{R}^{J}$. Then Corollary \ref{cor5.1} also holds for
$\tilde{W}$.
\end{remark}

\begin{proposition}
(viscosity subsolution) The function $\tilde{V}\left(  t,x,\hat{p},q\right)  $
is a viscosity subsolution of HJI equation \eqref{5.1}.
\end{proposition}

\paragraph{Proof}
For fixed $(\hat p,q)\in\mathbf{R}^I\times\Delta(J)$, we denote $\tilde V(t,x,\hat p,q)$  by $\tilde V(t,x)$ for short.
Fix $(t,x)\in[0,T]\times\mathbf{R}^d$. 
Let $M>0$ such that, for all $(s,y)\in[0,T]\times\overline B_1(x)$, all $n\in\mathbf{N}^*$ and all 
$(\alpha,\beta)\in{\cal A}^{\pi_n}_r\times{\cal B}^{\pi_n}_r$,
$X^{t,y,\alpha,\beta}_s\in \overline B_M(x)$ (where $\overline B_M(x)$ is the closed ball of center $x$ and radius $M$).\\
Let $\varphi\in{\cal C}^1([0,T]\times\mathbf{R}^d)$  such that $\tilde V-\varphi$ attains a strict maximum on $[0,T]\times\overline B_M(x)$ at $(t,x)$:
\[ (\tilde V-\varphi)(t,x)>(\tilde V-\varphi)(s,y),\; (s,y)\in[0,T]\times\overline B_M(x)\setminus\{ (t,x)\}.\]
For all $n\in\mathbf{N}$, let $(s_n,x_n)\in[0,T]\times\overline B_M(x)$ be such that $V^{\pi_n*}-\varphi$ achieves at $(s_n,x_n)$ its maximum over $[0,T]\times\overline B_M(x)$. 
Then there exists a subsequence (still denoted by $(s_n,x_n)$ which converges to $(t,x)$ (see \cite{cil} Proposition 4.3).\\
Let $N_0\in\mathbf{N}$ such that, for all $n\geq N_0$, $|x-x_n|\leq 1$. 
Then, for any $n\geq N_0$, let $k_n$ be such that $t^n_{k_{n-1}}\leq s_n< t^n_{k_n}$.
By the dynamic programming principle (Lemma 4.1), we have
\[ \begin{array}{rl}
0\leq & \inf_{\beta\in{\cal B}^{\pi_n}_r(s_n,t^n_{k_n})}\sup_{\alpha\in{\cal A}^{\pi_n}_r(s_n,t^n_{k_n})}E[V^{\pi_n*}(t^n_{k_n},X^{s_n,x_n,\alpha,\beta}_{t_{k_n}})-V^{\pi_n*}(s_n,x_n)]\\

\leq &\inf_{\beta\in{\cal B}^{\pi_n}_r(s_n,t^n_{k_n})}\sup_{\alpha\in{\cal A}^{\pi_n}_r(s_n,t^n_{k_n})}E[\varphi(t^n_{k_n},X^{s_n,x_n,\alpha,\beta}_{t_{k_n}})-\varphi(s_n,x_n)]\\

=&\inf_{\beta\in{\cal B}^{\pi_n}_r(s_n,t^n_{k_n})}\sup_{\alpha\in{\cal A}^{\pi_n}_r(s_n,t^n_{k_n})}
E\big[\int_{s_n}^{t^n_{k_n}}\big(\frac{\partial}{\partial r}\varphi(r,X^{s_n,x_n,\alpha,\beta}_r)\\
&\qquad\qquad\qquad\qquad
+f(X^{s_n,x_n,\alpha,\beta}_r,\alpha_r,\beta_r)D\varphi(r,X^{s_n,x_n,\alpha,\beta}_r)\big)dr\big]
\\
\leq &
\inf_{\beta\in{\cal B}^{\pi_n}_r(s_n,t^n_{k_n})}\sup_{\alpha\in{\cal A}^{\pi_n}_r(s_n,t^n_{k_n})}
E\big[\int_{s_n}^{t^n_{k_n}}\big(\frac{\partial}{\partial r}\varphi(s_n,x_n)\\
&\qquad\qquad\qquad\qquad
+f(x_n,\alpha_r,\beta_r)D\varphi(s_n,x_n)+m(C|t^n_{k_n}-s_n|)\big)dr\big],
\end{array}\]
with 
\[ \begin{array}{l}
m(\delta):=\\
\sup_{E_\delta}
\left|\left(\frac{\partial}{\partial r}\varphi(s,y)+f(y,u,v)D\varphi(s,y)\right)
-\left(\frac{\partial}{\partial r}\varphi(r,y')+f(y',u,v)D\varphi(r,y')\right)\right|,
\end{array}\]
with $E_\delta=\{u\in U,v\in V,s,r\in[0,T],y,y'\in \overline B_M(x); |s-r|+|y-y'|\leq \delta\}$.
Remark that $\lim_{\delta\searrow 0}m(\delta)=0$.\\
It follows that 
\[\begin{array}{l}
-(t^n_{k_n}-s_n)\left(\frac{\partial}{\partial t}\varphi(s_n,x_n)+m(C|t^n_{k_n}-s_n|)\right)\\
\qquad\qquad\leq 
\inf_{\beta\in{\cal B}^{\pi_n}_r(s_n,t^n_{k_n})}\sup_{\alpha\in{\cal A}^{\pi_n}_r(s_n,t^n_{k_n})}
E\left[\int_{s_n}^{t^n_{k_n}}f(x_n,\alpha_r,\beta_r)D\varphi(s_n,x_n)dr\right]\\
\qquad\qquad\leq \sup_{\alpha\in{\cal A}^{\pi_n}_r(s_n,t^n_{k_n})}
E\left[\int_{s_n}^{t^n_{k_n}}f(x_n,\alpha_r,\tilde\beta_r)D\varphi(s_n,x_n)dr\right],
\end{array}\]
where $\tilde\beta(\omega,u)=\sum_{l=k_n+1}^{N^n}\tilde v(\zeta_{2,l})(\omega)\mathbf{1}%
_{[t_{l-1},t_{l})}(s)$, for some arbitrary measurable map $\tilde v:[0,1]\rightarrow V$.\\
Now we can find a $(t^n_{k_n}-s_n)^2$-optimal strategy $\alpha^n$ (depending on $\tilde\beta$) such that
\[ \begin{array}{l}
-(t^n_{k_n}-s_n)\left(\frac{\partial}{\partial t}\varphi(s_n,x_n)+m(C|t^n_{k_n}-s_n|)+(t^n_{k_n}-s_n)\right)\\
\leq E\left[\int_{s_n}^{t^n_{k_n}}f(x_n,\alpha^n_r,\tilde\beta_r)D\varphi(t,x)dr\right]\\
=E\left[\int_{s_n}^{t^n_{k_n}}f(x_n,\alpha^n(\zeta_{1,k_n},\tilde v)_r,\tilde v(\zeta_{2,k_n}))D\varphi(s_n,x_n)dr\right].
\end{array}\]
Due to the time delay, on $[s_n,t^n_{k_n}]$, $\alpha^n$ doesn't depend on the control $\tilde v$ of Player 2.
 (Note that this is the crucial point where we use that both players have the same time grid.) Then, thanks to the independence between $\zeta_{1,k_n}$ and $\zeta_{2,k_n}$, we get
\[\begin{array}{l}
-(t^n_{k_n}-s_n)\left(\frac{\partial}{\partial t}\varphi(s_n,x_n)+m(C|t^n_{k_n}-s_n|)+(t^n_{k_n}-s_n)\right)\\
\leq (t^n_{k_n}-s_n)\sup_{\mu\in{\cal P}(U)}\int_UE[f(x_n,u,\tilde v(\zeta_{2,k_n}))D\varphi(s_n,x_n))\mu(du).
\end{array}\]
Thus, thanks to the arbitrariness of $\tilde v$,
\[\begin{array}{l}
-(t^n_{k_n}-s_n)\left(\frac{\partial}{\partial t}\varphi(s_n,x_n)+m(C|t^n_{k_n}-s_n|)+(t^n_{k_n}-s_n)\right)\\
\leq (t^n_{k_n}-s_n)\inf_{\tilde v}\sup_{\mu\in{\cal P}(U)}\int_UE[f(x_n,u,\tilde v(\zeta_{2,k_n}))D\varphi(s_n,x_n))\mu(du)\\
=(t^n_{k_n}-s_n)\inf_{\nu\in{\cal P}(V)}\sup_{\mu\in{\cal P}(U)}\int_{U\times V}f(x_n,u,v)D\varphi(s_n,x_n)\mu(du)\nu(dv).
\end{array}\]
Now recall that $(s_n,x_n)\rightarrow(t,x)$ and that $0\leq (t^n_{k_n}-s_n)\leq (t^n_{k_n}-t^n_{k_{n-1}})\leq |\pi_n|\rightarrow 0$ as $n\rightarrow\infty$.\\
Therefore, taking the limit as $n\rightarrow\infty$, we obtain
\[ -\frac{\partial}{\partial t}\varphi(t,x)\leq\inf_{\nu\in{\cal P}(V)}\sup_{\mu\in{\cal P}(U)}\int_{U\times V}f(x,u,v)D\varphi(t,x)\mu(du)\nu(dv).\]
The result follows.

 \begin{proposition}
(viscosity supersolution) The function $\tilde{W}\left(  t,x,\hat{p},q\right)  $
is a viscosity supersolution of HJI equation \eqref{5.1}.
\end{proposition}

\paragraph{Proof}
Note that
\[
-W^{\pi}\left(  t,x,p,q\right)  =\underset{\left(  \alpha_{i}\right)
\in\left(  \mathcal{A}_{r}^{\pi}\left(  t,T\right)  \right)  ^{I}}{\sup
}\underset{\left(  \beta_{j}\right)  \in\left(  \mathcal{B}_{r}^{\pi}\left(
t,T\right)  \right)  ^{J}}{\inf}\sum_{i=1}^{I}\sum_{j=1}^{J}p_{i}q_{j}E\left[
-g_{ij}\left(  X_{T}^{t,x,\alpha_{i},\beta_{j}}\right)  \right]  .
\]
The right side of this equation has the same form as $V^{\pi}$, only the role
of players changes. Hence, the convex conjugate of $\left(  -W^{\pi}\right)  $
with respect to $q$, i.e., $-W^{\pi\#}\left(  -\hat{q}\right)  $ satisfies a
subdynamic programming principle. Then, as a consequence of the above result
for $\tilde V$, we can deduce easily the following:\\
For any $\left(  t,x,p,\hat{q}\right)  \in
\lbrack0,T)\times\mathbf{R}^{d}\times\Delta\left(  I\right)  \times
\mathbf{R}^{J}$, and for all $l$ ($k\leq l\leq n$), we have%
\begin{equation}
W^{\pi_{n}\#}\left(  t,x,p,\hat{q}\right)  \geq\underset{\alpha\in
\mathcal{A}_{r}^{\pi}\left(  t,t_{k}^{n}\right)  }{\sup}\underset{\left(
\beta_{j}\right)  \in\mathcal{B}_{r}^{\pi}\left(  t,t_{k}^{n}\right)  }{\inf
}E\left[  W^{\pi_{n}\#}\left(  t_{k}^{n},X_{t_{k}^{n}}^{t,x,\alpha,\beta
},p,\hat{q}\right)  \right]  ,
\end{equation}
and $\tilde{W}$, the uniform limit on compact sets of $(W^{\pi_{n}\#})$, is a
supersolution of the HJI equation \eqref{5.1}.\\
(Here we have used equality \eqref{5.2}. Actually $\tilde{W}$ is a
supersolution of \eqref{5.1} with Hamiltonian%
\[
H^{\ast}\left(  x,\xi\right)  =\underset{\mu\in\mathcal{P}\left(  U\right)
}{\sup}\underset{\nu\in\mathcal{P}\left(  V\right)  }{\inf}%
{\displaystyle\int\nolimits_{U\times V}}
f\left(  x,u,v\right)  \mu\left(  du\right)  \nu\left(  dv\right)  \cdot\xi.)
\]
$\Box$

\section{Existence of the value}

In this section we show that the limit of the game along partitions has a
value. This value can be characterized by dual solutions of some HJI equation.

We now recall the definition of dual solutions for the following HJI
equation:
\begin{equation}
\left\{
\begin{array}
[c]{ll}%
\frac{\partial V}{\partial t}\left(  t,x\right)  +H\left(  x,DV\left(
t,x\right)  \right)  =0\text{,} & \text{in }\left[  0,T\right]  \times
\mathbf{R}^{d}\text{,}\\
V\left(  T,x\right)  =\sum_{ij}p_{i}q_{j}g_{ij}\left(  x\right)  \text{,} &
\end{array}
\right.  \label{6.1}%
\end{equation}
where
\[
H\left(  x,\xi\right)  =\underset{\mu\in\mathcal{P}\left(  U\right)  }{\inf
}\underset{\nu\in\mathcal{P}\left(  V\right)  }{\sup} \left(
{\displaystyle\int\nolimits_{U\times V}} f\left(  x,u,v\right)  \mu\left(
du\right)  \nu\left(  dv\right)  \cdot\xi\right)  .
\]

\begin{definition}
A function $w:[0,T]\times\mathbf{R}^{d}\times\Delta\left(  I\right)
\times\Delta\left(  J\right)  \mapsto\mathbf{R}$ is a dual subsolution of
\eqref{6.1} if w is Lipschitz continuous, convex with respect to p, and
concave with respect to q and if for any $\left(  p,\hat{q}\right)  \in
\Delta\left(  I\right)  \times$\textbf{$R$}$^{J}$, $\left(  t,x\right)
\mapsto w^{\#}\left(  t,x,p,\hat{q}\right)  $ is a supersolution of the dual
HJ equation%
\begin{equation}%
\begin{array}
[c]{ll}%
\frac{\partial V}{\partial t}\left(  t,x\right)  +H^{\ast}\left(  x,DV\left(
t,x\right)  \right)  =0\text{,} & \text{in }\left[  0,T\right]  \times
\mathbf{R}^{d}\text{,}%
\end{array}
\label{6.3}%
\end{equation}
where $H^{\ast}\left(  x,\xi\right)  =-H\left(  x,-\xi\right)  $.

We call $w:[0,T]\times\mathbf{R}^{d}\times\Delta\left(  I\right)  \times
\Delta\left(  J\right)  \mapsto\mathbf{R}$ a dual supersolution of
\eqref{6.1}, if w is Lipschitz continuous, convex with respect to p, and
concave with respect to q and if for any $\left(  \hat{p},q\right)
\in\mathbf{R}^{I}\times\Delta\left(  J\right)  $, $\left(  t,x\right)  \mapsto
w^{\ast}\left(  t,x,p,\hat{q}\right)  $ is a subsolution of \eqref{6.3}.

The function $w$ is called the dual solution of \eqref{6.1} if w is at the
same time a dual subsolution and a dual supersolution of \eqref{6.3}.
\end{definition}

Note that $\mathcal{P}\left(  U\right)  $ and $\mathcal{P}\left(  V\right)  $
are compact spaces. The measures $\mu\in\mathcal{P}\left(  U\right)  $ and
$\nu\in\mathcal{P}\left(  V\right)  $ can be interpreted as control variables.
Therefore, the comparison principle applies here in the sense of Cardaliaguet
\textrm{\cite{c07}}.

\begin{lemma}
\label{lem6.1}(comparison principle) Let $w_{1},w_{2}:\left[  0,T\right]
\times\mathbf{R}^{d}\times\Delta\left(  I\right)  \times\Delta\left(
J\right)  \mapsto\mathbf{R}$ a dual subsolution and a dual supersolution of HJ
\eqref{6.1}, respectively. If, for all $\left(  x,p,q\right)  \in
\mathbf{R}^{d}\times\Delta\left(  I\right)  \times\Delta\left(  J\right)  $,
$w_{1}\left(  T,x,p,q\right)  \leq w_{2}\left(  T,x,p,q\right)  $, then we
have $w_{1}\leq w_{2}$ on $\left[  0,T\right]  \times\mathbf{R}^{d}%
\times\Delta\left(  I\right)  \times\Delta\left(  J\right)  $.
\end{lemma}

We now state the main result of the paper.

\begin{theorem}
\label{unex} (uniqueness and existence of the value) For all sequences of
partitions $(\pi_{n})$ with $|\pi_{n}|\rightarrow0$, the sequences $\left(
V^{\pi_{n} }\right)  $ and $\left(  W^{\pi_{n}}\right)  $ converge uniformly
on compact sets to a same Lipschitz continuous function $V$, which is the
unique dual solution of the HJ equation \eqref{6.1}.
\end{theorem}

We will establish a, in appearance, weaker result :

\begin{proposition}
\label{thm6.1} For all sequences of partitions $(\pi_{n})$ with $|\pi
_{n}|\rightarrow0$, there exists a subsequence of partitions, still denoted by
$\left(  \pi_{n}\right)  _{n\geq1}$, such that $(V^{\pi_{n}},W^{\pi_{n}})$
converges uniformly on compact sets to a couple $(V,V)$, where the function
$V$ is the unique solution of the HJ equation \eqref{6.1}.
\end{proposition}

But we remark that, if Proposition \ref{thm6.1} is true for the partition
$(\pi_{n})$, it holds also for all subsequence of $(\pi_{n})$: there exists a
sub-subsequence $(\pi_{n_{l}})$ such that $(V^{\pi_{n_{l}}},W^{\pi_{n_{l}}})$
converges uniformly on compact sets. But Proposition \ref{thm6.1}
characterizes the limit $V(=W)$ as the unique dual solution of the
Hamilton-Jacobi-equation \eqref{6.3}. Consequently, all converging
sub-subsequences have the same limit, and the Theorem \ref{unex}
follows.\newline

\noindent\textbf{Proof of Proposition \ref{6.1}}\newline\indent Step 1. We
know since Lemma 3.1 that, for all $n$, $V^{\pi_{n}}$ and $W^{\pi_{n}}$ are
Lipschitz continuous in all their variables, with the same Lipschitz constant
depending only on $f$ and $g$. It follows that there is a subsequence of
partitions, still denoted by $\left(  \pi_{n}\right)  _{n\geq1}$ and two
bounded Lipschitz functions $V$, $W$:$\left[  0,T\right]  \times\mathbf{R}%
^{d}\times\Delta\left(  I\right)  \times\Delta\left(  J\right)  \mapsto
\mathbf{R}$ such that $\left(  V^{\pi_{n}},W^{\pi_{n}}\right)  \rightarrow
\left(  V,W\right)  $ uniformly on compact sets in $\left[  0,T\right]
\times\mathbf{R}^{d}\times\Delta\left(  I\right)  \times\Delta\left(
J\right)  $ and, obviously, the functions $V,W$ are also Lipschitz continuous
with respect to all their variables, convex in $p$, concave in $q$ (see Lemma
3.2).\newline

Step 2. Recall that $\tilde{W}=$ $\underset{n\rightarrow+\infty}{\lim}\left(
W^{\pi_{n}}\right)  ^{\#}$, $\tilde{V}=\underset{n\rightarrow+\infty}{\lim
}\left(  V^{\pi_{n}}\right)  ^{\ast}$. Due to our results in Section 5,
$w_{1}:=\tilde{W}^{\#}$ and $w_{2}=\tilde{V}^{\ast}$ are, respectively, a dual
subsolution and a dual supersolution of HJI equation \eqref{6.1} and
$\tilde{W}^{\#}\left(  T,x,p,q\right)  =\tilde{V}^{\ast}\left(
T,x,p,q\right)  =\sum_{ij}p_{i}q_{j}g_{ij}\left(  x\right)  $. Then, by the
comparison principle (Lemma \ref{lem6.1}), we have%
\begin{equation}
\tilde{W}^{\#}\leq\tilde{V}^{\ast}\text{ on }\left[  0,T\right]
\times\mathbf{R}^{d}\times\Delta\left(  I\right)  \times\Delta\left(
J\right)  . \label{6.5}%
\end{equation}

Step 3. Since $V^{\pi_{n}}$ converges to $V$ uniformly on compact sets, for
any $\varepsilon,M>0$, there is a positive integer $N_{\varepsilon,M}$ such
that for all $n\geq N_{\varepsilon,M}$, $t\in\left[  0,T\right]  $ and
$\left\vert x\right\vert \leq M$, $\left(  p,q\right)  \in\Delta
(I)\times\Delta\left(  J\right)  $,
\[
\left\vert V^{\pi_{n}}\left(  t,x,p,q\right)  -V\left(  t,x,p,q\right)
\right\vert \leq\varepsilon.
\]
Here $\varepsilon>0$ is arbitrarily given. Consequently, we have for all
$n\geq N_{\varepsilon,M}$, $t\in\left[  0,T\right]  $, $\left\vert
x\right\vert \leq M$ and $(\hat p,q)\in\mathbf{R}^{I}\times\Delta(J)$,
\begin{align*}
\left\vert V^{\pi_{n}*}\left(  t,x,\hat p,q\right)  -V^{\ast}\left(  t,x,\hat
p,q\right)  \right\vert  &  \leq\left\vert \underset{p\in\Delta(I)}{\sup
}\left\{  \hat{p}\cdot p-V^{\pi_{n}}\left(  t,x,\hat{p},q\right)  \right\}
-\underset{p\in\Delta(I)}{\sup}\left\{  \hat{p}\cdot p-V\left(
t,x,p,q\right)  \right\}  \right\vert \\
&  \leq\underset{p\in\Delta(I)}{\sup}\left\vert V^{\pi_{n} }\left(
t,x,\hat{p},q\right)  -V\left(  t,x,p,q\right)  \right\vert \\
&  \leq\varepsilon.
\end{align*}
Hence, $\tilde{V}=\underset{n\rightarrow+\infty}{\lim}V^{\pi_{n}*}=V^{*}$,
and, consequently, since $V$ is convex in $p$, $V=V^{**}=\tilde V^{*}$.

In a symmetric way, it is easy to get that $\tilde{W}^{\#}=\underset
{n\rightarrow+\infty}{\lim}W^{\pi_{n}}=W$.

As a consequence of \eqref{6.5}, we have $W\leq V$ in $\left[  0,T\right]
\times\mathbf{R}^{d}\times\Delta\left(  I\right)  \times\Delta\left(
J\right)  $. Observing that $W=$ $\underset{n\rightarrow+\infty}{\lim}%
W^{\pi_{n}}\geq$ $\underset{n\rightarrow+\infty}{\lim}V^{\pi_{n}}=V$, we
obtain that the game has a limit value%
\[
W=V\text{ in }\left[  0,T\right]  \times\mathbf{R}^{d}\times\Delta\left(
I\right)  \times\Delta\left(  J\right)  .
\]

From the above proof we deduce that the value $V$ $\left(  =W\right)  $ is the
unique dual solution of HJ equation \eqref{6.1}. $\Box$

\section{The case of lack of information on the dynamics}

In this section, we consider a game in which the players have also an
asymmetric information on the dynamics: In each of the scenarios $\left(
i,j\right)  \in\left\{  1,\ldots,I\right\}  \times\left\{  1,\ldots,J\right\}
$ the game has a different dynamic given by%
\begin{equation}
\left\{
\begin{array}
[c]{ll}%
dX_{ij}^{t,x_{ij},u,v}\left(  s\right)  =f_{ij}\left(  X_{ij}^{t,x_{ij}%
,u,v}\left(  s\right)  ,u\left(  s\right)  ,v\left(  s\right)  \right)
ds\text{,} & s\in\left[  t,T\right]  \text{, }\left(  u,v\right)
\in\mathcal{U}_{t,T}\times\mathcal{V}_{t,T}\\
X_{ij}^{t,x,u,v}\left(  t\right)  =x_{ij}\in\mathbf{R}^{d}\text{,} &
\end{array}
\right.
\end{equation}
where $f_{ij}:\mathbf{R}^{d}\times U\times V\mapsto\mathbf{R}^{d}$ is bounded,
continuous and Lipschitz continuous in $x$, uniformly with respect to $(u,v)$.
The possible payoffs are still given by $g_{ij}\left(  X_{ij}^{t,x_{ij}%
,u,v}\left(  T\right)  \right)  $, with $g_{ij}$ bounded and Lipschitz.

As already in the sections before, for some fixed $(p,q)\in\Delta
(I)\times\Delta(J)$, at time $t$, a scenario ($i,j$) is chosen at random with
the probability $p_{i}q_{j}$; the choice of $i$ is communicated only to player
1, while the choice of $j$ is communicated only to player 2. The players
observe their opponent's behavior and try to deduce from it their missing
information. Player 1 aims to minimize, Player 2 to maximize the payoff. We
remark that, in the case with Isaacs' assumption but with correlated
information, the case of lack of information on the dynamics has been
considered by Oliu Barton \cite{oliubarton}. \newline

As in the previous chapters, we fix a partition $\pi=\{ t_{0}=t<\ldots
<t_{N}=T\}$. Keeping the same definitions for the sets of strategies
$\mathcal{A}_{r}^{\pi}$ and $\mathcal{B}_{r}^{\pi}$, the upper and lower
values of the game associated to $\pi$ are then, respectively, for $\left(
p,q\right)  \in\Delta\left(  I\right)  \times\Delta\left(  J\right)  $,
$\left(  t,\mathbf{x}\right)  \in\lbrack0,T)\times\left(  \mathbf{R}%
^{d}\right)  ^{IJ}$,%

\begin{equation}
W^{\pi}\left(  t,\mathbf{x},p,q\right)  =\underset{\left(  \alpha_{i}\right)
\in\left(  \mathcal{A}_{r}^{\pi}\left(  t,T\right)  \right)  ^{I}}{\inf
}\underset{\left(  \beta_{j}\right)  \in\left(  \mathcal{B}_{r}^{\pi}\left(
t,T\right)  \right)  ^{J}}{\sup}\sum_{i=1}^{I}\sum_{j=1}^{J}p_{i}q_{j}E\left[
g_{ij}\left(  X_{ij}^{t,x_{ij},\alpha_{i},\beta_{j}}\left(  T\right)  \right)
\right]  \text{, }%
\end{equation}
and
\begin{equation}
V^{\pi}\left(  t,\mathbf{x},p,q\right)  =\underset{\left(  \beta_{j}\right)
\in\left(  \mathcal{B}_{r}^{\pi}\left(  t,T\right)  \right)  ^{J}}{\sup
}\underset{\left(  \alpha_{i}\right)  \in\left(  \mathcal{A}_{r}^{\pi}\left(
t,T\right)  \right)  ^{I}}{\inf}\sum_{i=1}^{I}\sum_{j=1}^{J}p_{i}q_{j}E\left[
g_{ij}\left(  X_{ij}^{t,x_{ij},\alpha_{i},\beta_{j}}\left(  T\right)  \right)
\right]  \text{. }%
\end{equation}

The idea to solve this case of asymmetric information is to blow up the
dynamics: We introduce the following auxiliary dynamic with values in
$\mathbf{R}^{d\times I\times J}$:
\begin{equation}
\left\{
\begin{array}
[c]{ll}%
d\mathbf{X}^{t,\mathbf{x},u,v}\left(  s\right)  =\mathbf{F}\left(
\mathbf{X}^{t,\mathbf{x},u,v}\left(  s\right)  ,u\left(  s\right)  ,v\left(
s\right)  \right)  ds\text{,} & s\in\left[  t,T\right]  \text{,}\\
\mathbf{X}^{t,\mathbf{x},u,v}\left(  t\right)  =\mathbf{x}\text{,} &
\end{array}
\right.
\end{equation}
with $\mathbf{F=}\left(  f_{ij}\right)  ,\mathbf{x=}\left(  x_{ij}\right)  $
and a new family of payoffs $G_{ij}\left(  \mathbf{x}\right)  =g_{ij}\left(
x_{ij}\right)  $. This permits us to rewrite the value functions as%

\[
W^{\pi}\left(  t,\mathbf{x},p,q\right)  =\underset{\left(  \alpha_{i}\right)
\in\left(  \mathcal{A}_{r}^{\pi}\left(  t,T\right)  \right)  ^{I}}{\inf
}\underset{\left(  \beta_{j}\right)  \in\left(  \mathcal{B}_{r}^{\pi}\left(
t,T\right)  \right)  ^{J}}{\sup}\sum_{i=1}^{I}\sum_{j=1}^{J}p_{i}q_{j}E\left[
G_{ij}\left(  \mathbf{X}_{T}^{t,\mathbf{x},\alpha_{i},\beta_{j}}\right)
\right]
\]
and
\[
V^{\pi}\left(  t,\mathbf{x},p,q\right)  =\underset{\left(  \beta_{j}\right)
\in\left(  \mathcal{B}_{r}^{\pi}\left(  t,T\right)  \right)  ^{J}}{\sup
}\underset{\left(  \alpha_{i}\right)  \in\left(  \mathcal{A}_{r}^{\pi}\left(
t,T\right)  \right)  ^{I}}{\inf}\sum_{i=1}^{I}\sum_{j=1}^{J}p_{i}q_{j}E\left[
G_{ij}\left(  \mathbf{X}_{T}^{t,\mathbf{x},\alpha_{i},\beta_{j}}\right)
\right]  ,\newline%
\]
and we recover the case of asymmetric information solved in the previous
chapters. As a consequence we have the following result:

\begin{theorem}
\label{thm6.2} The limit value of the game exists, as the mesh of partitions
tends zero, and it is the dual solution of the following HJI equation:
\begin{equation}
\left\{
\begin{array}
[c]{ll}%
\frac{\partial V}{\partial t}\left(  t,\mathbf{x}\right)  +\mathbf{H}\left(
\mathbf{x},DV\left(  t,\mathbf{x}\right)  \right)  =0\text{,} & \text{in
}\left[  0,T\right]  \times\left(  \mathbf{R}^{d}\right)  ^{IJ}\text{,}\\
V\left(  T,\mathbf{x}\right)  =\sum_{ij}p_{i}q_{j}g_{ij}\left(  x_{ij}\right)
\text{,} &
\end{array}
\right.
\end{equation}
where
\begin{align*}
\mathbf{H}\left(  \mathbf{x},\xi\right)   &  =\underset{\mu\in\mathcal{P}%
\left(  U\right)  }{\inf}\underset{\nu\in\mathcal{P}\left(  V\right)  }{\sup
}\sum_{ij}%
{\displaystyle\int\nolimits_{U\times V}}
f_{ij}\left(  \mathbf{x}_{ij},u,v\right)  \mu\left(  du\right)  \nu\left(
dv\right)  \cdot\xi_{ij}\\
&  =\underset{\nu\in\mathcal{P}\left(  V\right)  }{\sup}\underset{\mu
\in\mathcal{P}\left(  U\right)  }{\inf}\sum_{ij}%
{\displaystyle\int\nolimits_{U\times V}}
f_{ij}\left(  \mathbf{x}_{ij},u,v\right)  \mu\left(  du\right)  \nu\left(
dv\right)  \cdot\xi_{ij}.
\end{align*}

\end{theorem}

\begin{remark}
In the present paper the asymmetric information structure concerns a number of
finite types of the payoffs $g_{ij}(X_{T})$, $i=1,2\ldots I$, $j=1,2\ldots J$.
An interesting question would concern the case where the set of types is
infinite : typically each player, instead knowing a probability measure with
finite support which represent the information on its opponents knowledge,
would know a probability measure which support is an arbitrary subset of some
$\mathbf{R}^{q}$.

In this direction the results of the present paper could be extended using
methods of \cite{cr12} and \cite{CJQ}.
\end{remark}

 \section*{Acknowledgments.}
 We thank the anonymous referees for the careful reading of the manuscript and the useful remarks and comments.

\end{document}